\documentclass[a4paper,12pt]{amsproc}
\usepackage{amsmath,amsfonts,amssymb,amsthm,anysize}
\usepackage{enumerate}
\usepackage{epsfig}
\usepackage{supertabular}
\usepackage{graphicx}

\newcommand\N{{\mathbb N}}
\newcommand\Z{{\mathbb Z}}
\newcommand\Q{{\mathbb Q}}

\newcommand\C{{\mathbb C}}
\newcommand\F{{\mathbb F}}

\newcommand\Gal{{\mathrm{Gal}}}

\newcommand\Tr{{\mathrm{Tr}}}
\newcommand\Norm{{\mathrm{Norm}}}

\newcommand\ord{{\mathrm{ord}}}

\newcommand\lcm{{\mathrm{lcm}}}

\newcommand\ann{\mathsf{ann}}

\renewcommand\mod{{\mathrm{mod\, \, }}}
\newcommand\pro{{\bf Proof: }}

\newcommand\D{\mathsf{D}}

\newcommand\OO{\mathsf{D}}

\newcommand\XX{\mathsf{X}}
\newcommand\DD{\mathsf{D}}
\newcommand\AAA{\mathsf{A}}

\theoremstyle{plain}
\newtheorem{theorem}{Theorem}[section]
\newtheorem{problem}[theorem]{Problem}
\newtheorem{lemma}[theorem]{Lemma}
\newtheorem{corollary}[theorem]{Corollary}

\newtheorem{proposition}[theorem]{Proposition}

\numberwithin{equation}{section}
\allowdisplaybreaks[4]

\theoremstyle{remark}
\newtheorem{remark}[theorem]{Remark}

\usepackage{url, hyperref}
\hypersetup{citecolor=blue, linkcolor=blue, colorlinks=true}

\renewcommand\le{\leqslant}
\renewcommand\ge{\geqslant}

\begin{document}

\title[Powers of Gauss sums in quadratic fields]{Powers of Gauss sums in quadratic fields}


%
\author{Koji Momihara}
\address{ %
Division of Natural Science\\
Faculty of Advanced Science and Technology\\
Kumamoto University\\
2-40-1 Kurokami, Kumamoto 860-8555, Japan}
\email{momihara@educ.kumamoto-u.ac.jp}
\thanks{The author acknowledges the support by 
JSPS under Grant-in-Aid for Scientific Research (C) 20K03719.}


\begin{abstract}
In the past two decades, many researchers have 
studied {\it index $2$} Gauss sums, where the group generated by the characteristic $p$ of the underling finite field is of index $2$ in the unit group of $\Z/m\Z$ for the order $m$ of the multiplicative character involved.  
A complete solution to the problem of evaluating index $2$ Gauss sums was given by Yang and Xia~(2010). In particular, it is known that some nonzero integral powers of the Gauss sums in this case are in quadratic fields.  
On the other hand, Chowla~(1962), McEliece~(1974), Evans~(1977, 1981) and Aoki~(1997, 2004, 2012)  studied {\it pure} Gauss sums, some nonzero integral powers of which are in  the field of rational numbers. 

In this paper, we study Gauss sums, some integral powers of which are in quadratic fields. 
This class of Gauss sums is a
generalization of index $2$ Gauss sums and an extension of pure Gauss sums to quadratic fields. 
%
\end{abstract}


\maketitle

\section{Introduction}\label{sec:intro}
Let $p$ be a prime, $h$ be a positive integer, and $q=p^h$. The canonical additive character $\psi$ of $\F_q$ is defined by 
$$\psi\colon\F_q\to \C^{\ast},\qquad\psi(x)=\zeta_p^{\Tr _{q/p}(x)},$$
where $\zeta_p={\rm exp}(\frac {2\pi i}{p})$ and $\Tr _{q/p}$ is the trace from $\F_q$ to $\F_p$. For a multiplicative character  
$\eta_m$ of order $m$ of $\F_q$, we define the {\it Gauss sum} 
\[
G_{q}(\eta_m)=\sum_{x\in \F_q^\ast}\eta_m(x)\psi(x), 
\] 
which belongs to the ring of integers in the field $\Q(\zeta_m,\zeta_{p})$.

The Gauss sum is one of important and fundamental objects in number
theory. 
The concept of the Gauss sums was introduced by Gauss in 1801~\cite{Gauss}, who  evaluated the {\it quadratic} Gauss sums. 
\begin{theorem}{\em (\cite{Gauss})} \label{le:quad}
Let $\eta$ be the quadratic character of $\F_q=\F_{p^h}$. Then, it holds that 
\[
G_{q}(\eta)=(-1)^{h-1}\left(\sqrt{(-1)^{\frac{p-1}{2}}p}\right)^h. 
\]
\end{theorem}
After Gauss' work, many researchers have tried to evaluate Gauss sums for larger $m$. 
However, in general, the explicit evaluation of Gauss sums is a very difficult problem. 
There are only a few cases where the Gauss sums have been completely evaluated. 
For example,  the Gauss sums for $m = 3, 4, 5, 6, 8, 12,16,24$ have been 
evaluated (but not explicit in some cases). See \cite{BEW97} for more details. 

The next important case is the so-called
{\it semi-primitive case} (also referred to as uniform cyclotomy or supersingular), where there
exists an integer $s$ such that $p^s \equiv -1 \,(\mod{m})$. 
\begin{theorem}\label{thm:semiprim}{\em (\cite{BEW97})} 
Suppose that $m>2$ and $p$ is semi-primitive modulo $m$, 
i.e., there exists an $s$ such that  $p^s\equiv -1\,(\mod{m})$. Choose 
$s$ minimal and write 
$h=2st$. Let $\eta_m$ be a multiplicative character of order $m$. 
Then, 
\[
p^{-h/2}G_{p^h}(\eta_m)=
\left\{
\begin{array}{ll}
(-1)^{t-1}&  \mbox{if $p=2$;}\\
(-1)^{t-1+(p^s+1)t/m}&  \mbox{if $p>2$. }
 \end{array}
\right.
\]
\end{theorem}
The next interesting case is the {\it index $2$ case}, where the subgroup $\langle p\rangle$ generated by $p\in \Z$
has index $2$ in $(\Z/m\Z)^\times$. 
In this case, it is known that $m$ can have at most two distinct odd
prime divisors. Furthermore, some nonzero integral powers of index $2$ Gauss sums are in quadratic fields. Many authors have studied  this case, see, e.g., \cite{Lang,M98,M03,YX09,YX10}.
In particular, a complete solution to the problem of evaluating Gauss sums in this case was given in \cite{YX10}. As a large generalization, Aoki~\cite{A10} studied  Gauss sums such that $(\Z/m\Z)^\times /\langle p\rangle$ is  an  elementary abelian $2$-group. The {\it index $4$ case} including the case where $(\Z/m\Z)^\times /\langle p\rangle$  is cyclic was also studied in \cite{FY,FYL,YLF}.  

On the other hand, there were studies on Gauss sums from another point of view. Chowla~\cite{C01,C02} showed that if a Gauss sum defined in a prime field has the form $\epsilon p^{\frac{1}{2}}$ with $\epsilon$ a root of unity, it is in the quadratic case. 
McEliece~\cite{Mc} studied for which $(m,p,h)$, some nonzero integral power of the corresponding Gauss sum is an integer, i.e., $p^{-h/2}G_{p^h}(\eta_m)$ is a root of unity, related to weight distribution of irreducible cyclic codes.  Such Gauss sums are called {\it pure}. 
It is clear that the semi-primitive Gauss sums are examples of pure Gauss sums. 
Evans~\cite{E77} showed that  pure Gauss sums for prime powers $m$ are in the semi-primitive case. 
Furthermore,  he~\cite{E81} gave the following nontrivial families of pure Gauss sums which are not semi-primitive. 
\begin{theorem}\label{thm:small}
Suppose that $m=cd$ with $\gcd{(c,d)}=\gcd{(\ord_{c}(p),\ord_{d}(p))}=1$ and let $f=\ord_{m}(p)$, where $\ord_{n}(x)$ is the order of $x$ in $(\Z/n\Z)^\times$. Then, $G_{p^f}(\eta_m)$ is pure if any of the following holds. 
\begin{itemize}
\item[(1)] $\ord_{c}(p)=\phi(c)$ and $\ell\in \langle p\rangle\,(\mod{d})$ for some prime $\ell\,|\,c$. 
\item[(2)] $-1\not\in \langle p\rangle\,(\mod{c})$, $2\ord_{c}(p)=\phi(c)$, $\ell\in \langle p\rangle\,(\mod{d})$ for some prime $\ell\,|\,c$, and all of them hold with $c$ and $d$ interchanged. 
\item[(3)] $2||m$, $2+m/2\not\in \langle p\rangle\,(\mod{c})$,  $2\ord_{c}(p)=\phi(c)$, $-1$ or $\ell$ is in $\langle p\rangle\,(\mod{d})$ for some prime $\ell\,|\,c$, and all of them hold with $c$ and $d$ interchanged. 
\end{itemize} 
Here, $\phi$ is Euler's totient function.  
\end{theorem}
In particular, Aoki~\cite[Theorem~7.2]{A04} proved that the converse of the assertion of Theorem~\ref{thm:small} also  holds if $c$ and $d$ are both odd prime powers. 

As a remarkable work, Aoki~\cite{A04,A12} gave the following 
necessary and sufficient condition for a Gauss sum to be pure. 
\begin{theorem}{\em (\cite[Proposition~2.1]{A04})}\label{thm:aokip}
Let  $f=\ord_{m}(p)$. Then, 
 $G_{p^f}(\eta_m)$ is pure if and only if 
$\prod_{\ell\,|\,m}(1-\chi(\ell))=0$ for any odd Dirichlet character $\chi$ 
modulo $m$ such that $\chi(p)=1$, where $\ell$ runs over the prime divisors of $m$ not 
dividing the conductor of $\chi$. 
\end{theorem}
Based on this result, Aoki proved the following theorem on the finiteness of pure Gauss sums. 
\begin{theorem}{\em (\cite[Theorem~1.2]{A04})}\label{thm:aokifinite}
For a fixed $f$, the set of pairs $(m,[p]_m)$ such that $\ord_{m}(p)=f$ and  $G_{p^f}(\eta_m)$ is pure but not semi-primitive is finite, where $[p]_m$ is an integer 
such that $[p]_m\equiv p\,(\mod{m})$ and  $1\le [p]_m\le m-1$. 
\end{theorem}
Furthermore, he completely classified such pairs $(m,[p]_m)$ for $f=1,2,3,4$~\cite{A97}. This result was extended to $f=5,7,9,11,13,17,19,23$ in \cite{Mo}. 

As mentioned above, previous studies on evaluating Gauss sums fall into the following   
cases. 
\begin{itemize}
\item[(1)] $m=2,3,4,5,6,8,12,16,24$. 
\item[(2)] $\phi(m)/f=2,4$ or $(\Z/m\Z)^\times /\langle p\rangle$ is an elementary 
abelian $2$-group.  
\item[(3)] $G_{p^f}(\eta_m)/p^{f/2}$ is a root of unity. 
\end{itemize}
In this paper, we consider a new class of Gauss sums, some nonzero integral 
powers of which are in quadratic fields.  
This class of Gauss sums is a generalization of index $2$ Gauss sums and 
an extension of pure Gauss sums to quadratic fields. 
Thus, this class of Gauss sums unifies them.  
In particular, we are interested in classifying triples $(m,p,f)$ such that 
$G_{p^f}(\eta_m)^\ell$ is in a quadratic field for some nonzero integer $\ell$ 
but not concerned with explicitly evaluating the corresponding  Gauss sums. 

The objectives of this paper are three-fold.  
First, we give a general necessary  condition for some nonzero integral 
power of a Gauss sum to be in an extension field of minimal degree $e$ over $\Q$ for a fixed $e\ge 1$ in terms of Dirichlet characters modulo $m$. 
See Lemma~\ref{prop:equiv1} and Proposition~\ref{thm:mainpre} in Section~\ref{sec:nece}. The approach is based on the Stickelberger theorem on ideal factorization of 
Gauss sums and evaluating the exponent of each prime ideal factor of the Gauss sum by 
using $1$st generalized Bernoulli numbers, which was already applied by Aoki for pure Gauss sums. In this paper, we generalize this approach to general $e$. In particular, we show that the necessary condition is also sufficient if $e=2$. 
The following is one of our main theorems.  
\begin{theorem}
\label{cor:equiv2}
Let $f=\ord_{m}(p)$. 
Some nonzero integral power of 
$G_{p^f}(\eta_m)$ is in a quadratic field but any power of which is not in $\Q$
if and only if there exists a subgroup $E_0$ of index $2$ of $(\Z/m\Z)^\times$ containing $\langle p\rangle $ such that 
\begin{itemize}
\item[(1)] the unique nontrivial annihilator $\chi_{\ann}$ (as a Dirichlet character modulo $m$) of $E_0$ is an odd character and  $\chi_{\ann}(\ell)=-1$ for any prime divisor $\ell$ of  $m$ not dividing the conductor of $\chi_{\ann}$, and 
\item[(2)] 
$\prod_{\ell | m}(1-\chi(\ell))=0$ for any odd character $\chi$ modulo $m$ such that $\chi(p)=1$ except for $\chi_{\ann}$, where $\ell$ runs over all prime divisors of  $m$ not dividing the conductor of $\chi$. 
\end{itemize}
\end{theorem}
Second, we prove the  
following theorem in Section~\ref{sec:finite}, 
which is the quadratic field version of Theorem~\ref{thm:aokifinite}. 
\begin{theorem}\label{thm:aokifinite2}
For each fixed $f$, the set of pairs $(m,[p]_m)$ such that 
$\ord_{m}(p)=f$ and some nonzero integral power of 
$G_{p^f}(\eta_m)$ is in a quadratic field but any power of which is not in $\Q$ is finite. 
\end{theorem}
To prove Theorem~\ref{thm:aokifinite}, 
Aoki~\cite{A04} studied for which $m$ and $p$, the set of 
odd primitive characters modulo $m$ being trivial on $\langle p \rangle$ 
is empty in relation to Diophantine equations of  the form   $\sum_{i=1}^nx_i/d_i\equiv 0\,(\mod{1})$. 
On the other hand, we give a 
short proof for 
Theorem~\ref{thm:aokifinite2} 
based on the pigeonhole principle. 
Furthermore, we give a relatively sharp upper bound for $m$ in the case where $f$ is odd, and  completely classify such pairs $(m,[p]_m)$ for  $f=1,3,5,7$. 
However, we do not know whether this approach can be applied to general extension degree $e>2$. 

Third, we give a complete characterization for the new class of Gauss sums 
when 
$m$ is a prime power or $m$ has  exactly two distinct prime divisors with $8\not| m$ corresponding to Theorem~\ref{thm:small} and the converse of its assertion shown by Aoki~\cite{A04}. Then, we find infinitely many nontrivial examples of such Gauss sums not belonging to the index $2$ case.   

At last of this section, we mention again that we will not try to explicitly evaluate Gauss sums in this paper. To evaluate the Gauss sums in our case, one may need to deal with many cases individually as in the index $2$ case. Furthermore, we need to determine the minimum $\ell$ such that 
$G_{p^f}(\eta_m)^\ell\in F$ for some quadratic field $F$, which seems to be a difficult problem. Hence, this problem will be left as an open problem. (See 
Problem~\ref{prob:last} in the last section.) 
\section{Preliminaries}\label{sec:pre}
Let $p$ be a prime and let $q=p^h$ for a positive integer $h$. Let $m>1$ be a positive integer such that $m\,|\,q-1$ and  $\eta_m$ be a multiplicative character of order $m$ of $\F_q$.  

Let $\sigma_{a,b}$ be the automorphism of $\Q(\zeta_{m},\zeta_{p})$ determined 
by 
\[
\sigma_{a,b}(\zeta_m)=\zeta_{m}^a, \qquad
\sigma_{a,b}(\zeta_p)=\zeta_{p}^b 
\]
for $\gcd{(a,m)}=\gcd{(b,p)}=1$. 
Below are basic properties of Gauss sums \cite{LN97}: 
\begin{itemize}
\item[(i)] $G_q(\eta_m)\overline{G_q(\eta_m)}=q$ if $\eta_m$ is nontrivial;
\item[(ii)] $G_q(\eta_m^p)=G_q(\eta_m)$; 
\item[(iii)] $G_q(\eta_m^{-1})=\eta_m(-1)\overline{G_q(\eta_m)}$;
\item[(iv)] $G_q(\eta_m)=-1$ if $\eta_m$ is trivial;
\item[(v)] $\sigma_{a,b}(G_q(\eta_m))=\eta_m^{-a}(b)G_q(\eta_m^a)$. 
\end{itemize}
We give an important formula, known as 
 {\it the Davenport-Hasse lifting formula},  on Gauss sums.
\begin{theorem}\label{thm:lift}{\em (\cite[Theorem~11.5.2]{BEW97})}
Let $\eta_m$ be a nontrivial multiplicative character of $\F_q$ and 
let $\eta_m'$ be the lifted character of $\eta_m$ to $\F_{q^s}$, i.e., $\eta_m'(\alpha):=\eta_m(\Norm_{\F_{q^s}/\F_q}(\alpha))$ for $\alpha\in \F_{q^s}$. 
Then, it holds that 
\[
G_{q^s}(\eta_m')=(-1)^{s-1}(G_q(\eta_m))^s. 
\]
\end{theorem}

Hereafter, fix an integer $m>1$, and let  $p$ be a prime such that 
$\gcd{(p,m)}=1$. 
We will identify the quotient ring $\Z/m\Z$ with the ring of integers modulo $m$. Let $G=(\Z/m\Z)^\times$ be the unit group of $\Z/m\Z$, and 
let $f$ be the order of $p$ in $G$. Let $q=p^f$.

Define 
\[
K=\Q(\zeta_m),\, M=K(\zeta_p)=\Q(\zeta_m,\zeta_p), 
\]
and 
let $O_K$ and $O_M$ denote their respective rings of integers. 
For $j\in G$, define $\sigma_j\in \Gal(M/\Q(\zeta_p))$ by 
$\sigma_j(\zeta_m)=\zeta_m^j$. Let $P$ be a prime ideal of $O_K$ lying over 
$p$. Then, there is a prime ideal ${\mathfrak{p}}$ of $O_M$ such that 
$PO_M={\mathfrak{p}}^{p-1}$ and ${\mathfrak{p}}\cap O_K=P$. 
Denote $P_j=\sigma_j(P)$ and ${\mathfrak{p}}_j=\sigma_j({\mathfrak{p}})$,  and 
then $P_jO_M={\mathfrak{p}}_j^{p-1}$. Let $T$ be a set of representatives 
of
$G/\langle p\rangle$. Then, $pO_K=\prod_{j\in T}P_j$, 
where 
$P_j$ are all distinct. Hence, $pO_M=\prod_{j\in T}{\mathfrak{p}}_j^{p-1}$ holds. 

Define the character $\eta_P$ of order $m$ on the finite field 
$O_K/P$ by letting $\eta_P(\alpha+P)$ denote the 
unique power of $\zeta_m$ such that 
\[
\eta_P(\alpha+P)\equiv \alpha^{(q-1)/m}\, (\mod{P}), 
\] 
when $\alpha\in O_K\setminus P$. When $\alpha\in P$, let  $\eta_P(\alpha+P)=0$. 
Now we identify $\eta_P$ with a multiplicative character of 
$\F_q$. 

\begin{theorem}{\em (\cite[Theorems~11.2.2 and 11.2.7]{BEW97})} \label{thm:Stick}
 For a prime ideal  
${\mathfrak{p}}$ of $O_M$ lying over $P$, it holds that 
\[
G_q(\eta_P^{-1}) O_M=
{\mathfrak{p}}^{\frac{p-1}{m}\sum_{t\in T}\sum_{i=0}^{f-1}
[tp^i]_m\sigma_t^{-1}}.  
\] 
\end{theorem}
This theorem is known as {\it the Stickelberger theorem} on ideal factorizations of Gauss sums. 
On the other hand, the following is also known. 
\begin{theorem}{\em (\cite[Theorem~5.27]{LN97})} \label{thm:power}
It holds that 
$G_q(\eta_m)^m\in \Q(\zeta_m)$. 
\end{theorem}
Then, by the Stickelberger theorem and Theorem~\ref{thm:power}, we have 
\begin{equation}\label{eq:factor1}
G_q(\eta_P^{-1})^mO_K=P^{\sum_{t\in T}\sum_{i=0}^{f-1}
[tp^i]_m\sigma_t^{-1}}.  
\end{equation}
\section{Necessary and sufficient conditions}\label{sec:nece}
In this section, we give a necessary and sufficient condition for some nonzero integral  power of a Gauss sum to be in some quadratic field.
To do this, we define 
\begin{equation}\label{eq:defmu}
\mu_{p,f,m}:=\min\{[F:\Q]\,|\, \mbox{$G_{p^f}(\eta_P^{-1})^\ell\in F$ for some  $\ell\in \N$ and some field $F/\Q$}\}. 
\end{equation}
\begin{lemma}{\em (\cite{A04,E77,KL79})}\label{lem:pure1}
$\mu_{p,f,m}=1$ if and only if $\sum_{j=0}^{f-1}[tp^j]_m=fm/2$ for any integer $t$ prime to $m$. 
\end{lemma}
We give a general necessary condition below.  
\begin{lemma}\label{prop:equiv1}
There is a subgroup $E_0$ of index $\mu_{p,f,m}$ of $G$ containing $H=\langle p\rangle$ such that for each $i=0,1,\ldots,\mu_{p,f,m}-1$, $\sum_{j=0}^{f-1}[tp^j]_m$ is constant for any $t\in E_i$. 
Here,  
$R=\{s_i\,|\,i=0,1,\ldots,e-1\}$ is the set of coset representatives of 
$G/E_0$ and $E_i:=s_i E_0$ for $i=0,1,\ldots,\mu_{p,f,m}-1$. 
\end{lemma}
\pro 
We can take $\ell$ so that $m\,|\,\ell$. Then, $G_{p^f}(\eta_P^{-1})^\ell\in K$. This implies that $F$ is a  subfield of $K$. Furthermore, 
since $G_{p^f}(\eta_P^{-1})^\ell$ is invariant under $\sigma_{p}\in \Gal(K/\Q)$, $F$ is fixed by $\sigma_p$. 
Then, by \eqref{eq:factor1}, we have
\begin{align*}
P^{\frac{\ell}{m}\sum_{t\in G/H}\sum_{i=0}^{f-1}[tp^j]_m\sigma_t^{-1}}=&\,
(G_{p^f}(\eta_P^{-1})^\ell) O_{K}\\
=
&\,
\tau(G_{p^f}(\eta_P^{-1})^\ell)O_{K}=P^{\frac{\ell}{m}\sum_{t\in G/H}\sum_{i=0}^{f-1}[tp^j]_m\tau\sigma_t^{-1}}
\end{align*}
for any $\tau\in \Gal(K/F)$. Then, by comparing the exponents of $P^{\sigma_t}$, $t\in G$, of both sides, there is a subgroup $E_0$ of index  $\mu_{p,f,m}$ of $G$ such that for each $i=0,1,\ldots,\mu_{p,f,m}-1$, $\sum_{j=0}^{f-1}[tp^j]_m$ is constant for any $t\in E_i$ due to the Galois theory. 
In particular, $E_0$ contains $H=\langle p\rangle$ since $F$ is fixed by $\sigma_p$. 
\qed
\vspace{0.3cm}

Let $\DD(m)$ be the set of Dirichlet characters modulo $m$. Then, $\DD(m)$ can be viewed as the character group of $G$. For nontrivial $\chi\in \DD(m)$, we define the $1$st generalized Bernoulli number as 
$B_{1,\chi}:=\frac{1}{c_\chi}\sum_{x=0}^{c_{\chi}-1}x\chi(x)$, where $c_\chi$ is the conductor of $\chi$.  We use the following fact, see, e.g.,  \cite[Lemma~4.3]{A10}. 
\begin{lemma}\label{lem:ber}
For any nontrivial $\chi\in \DD(m)$, it holds that 
\[
\frac{1}{m}\sum_{x\in G}x\chi(x)=B_{1,\chi}\prod_{\ell | m}(1-\chi(\ell)), 
\]
where $\ell$ runs over all prime divisors of  $m$ not dividing the conductor of $\chi$.  
\end{lemma}

\begin{lemma}\label{lem:const1}
Let $E_0$ be a subgroup of index $e$ of $G$ containing $\langle p\rangle$. 
For each $i=0,1,\ldots,e-1$, 
if $\sum_{j=0}^{f-1}[tp^j]_m$ is constant, say $A_i$, for any $t\in E_i$, then $A_i=\frac{ef}{\phi(m)}\sum_{t\in E_i}t$. 
\end{lemma}
\pro
The condition $\sum_{j=0}^{f-1}[tp^j]_m=A_i$ implies that 
\[
|E_i|A_i=\sum_{t\in E_i}\sum_{j=0}^{f-1}[tp^j]_m=f\sum_{t\in E_i}t. 
\]
This completes the proof. 
\qed
\vspace{0.3cm}

Let  $\AAA(E_0)$ be the set of annihilators of $E_0$ in $\DD(m)$. The following is a generalization of \cite[Proposition~3.1]{A04}.  
\begin{proposition}\label{thm:mainpre}
Let $E_0$ be a subgroup of index $e$ of $G$ containing $H=\langle p\rangle$. 
Let $A_i=\frac{ef}{\phi(m)}\sum_{t\in E_i}t$, $i=0,1,\ldots, e-1$. 
For each $i$, $\sum_{j=0}^{f-1}[tp^j]_m=A_i$ for any $t\in E_i$ if and only if 
$\prod_{\ell | m}(1-\chi(\ell))=0$ for any odd character $\chi\in \DD(m)\setminus \AAA(E_0)$ such that $\chi(p)=1$, 
where $\ell$ runs over all prime divisors of $m$ not dividing the conductor of $\chi$.  
\end{proposition}
\pro 
Suppose that 
$\sum_{j=0}^{f-1}[tp^j]_m=A_i$ for any $t\in E_i$. 
Then, we have 
\begin{equation}\label{eq:tra1}
\sum_{t\in E_i}\left(\sum_{j=0}^{f-1}[tp^j]_m-A_i\right )\chi(t)=0, \quad \forall \chi \in \DD(m)\setminus \AAA(E_0). 
\end{equation}
Conversely, assume that \eqref{eq:tra1} holds. Then, 
for any $s\in E_i$, 
\begin{align}
0=&\sum_{\chi\in \DD(m)\setminus \AAA(E_0)}\sum_{t\in E_i}\left(\sum_{j=0}^{f-1}[tp^j]_m-A_i\right )\chi(t)\chi(s^{-1})\nonumber\\ 
=&\, \sum_{t\in E_i}\left(\sum_{j=0}^{f-1} [t p^j]_m-A_i\right )\sum_{\chi\in \DD(m)\setminus \AAA(E_0)}\chi(ts^{-1}).\label{eq:aaaa} 
\end{align}
Since 
\[
\sum_{\chi\in \DD(m)\setminus \AAA(E_0)}\chi(ts^{-1})=
\begin{cases}
\phi(m)-e&  \mbox{ if $t=s$,}\\
-e& \mbox{ otherwise, }
\end{cases}
\]
and  
$\sum_{t\in E_i}\left(\sum_{j=0}^{f-1} [t p^j]_m-A_i\right)=0$ as in the proof of Lemma~\ref{lem:const1}, 
we have 
\begin{align*}
\eqref{eq:aaaa}=&\, \phi(m)\left(\sum_{j=0}^{f-1}[s p^j]_m-A_i\right ).  
\end{align*}
Hence,  $\sum_{j=0}^{f-1}[sp^j]_m=A_i$ for any $s\in E_i$. 

We reformulate the left-hand side of \eqref{eq:tra1} as follows:  
\begin{align}
\sum_{t\in E_i}\left(\sum_{j=0}^{f-1}[tp^j]_m-A_i\right )\chi(t)\nonumber
=&\sum_{j=0}^{f-1}\sum_{t\in E_i}[tp^j]_m\chi(t)- A_i\sum_{t\in E_i}\chi(t)
\nonumber\\
=&\left(\sum_{j=0}^{f-1}\chi(p^{-j})\right)
\left(\sum_{t\in E_i}t\chi(t)\right)-A_i\sum_{t\in E_i}\chi(t).  \label{eq:tra2}
\end{align}
Note that the indicator function of each $E_i$ is given by 
\begin{equation}\label{eq:chara}
g_i(x)=\frac{1}{e}\sum_{\chi\in \AAA(E_0)}\chi(s_i^{-1}x). 
\end{equation}
Then, by \eqref{eq:chara} and Lemma~\ref{lem:ber}, 
we have
\begin{align*}
\eqref{eq:tra2}=&\,
\left(\sum_{j=0}^{f-1}\chi(p^{-j})\right)
\sum_{t\in G}t\chi(t)g_i(t) \\
=&\,\frac{m}{e}\left(\sum_{j=0}^{f-1}\chi(p^{-j})\right)
\sum_{\chi'\in \AAA(E_0)}{\chi'}^{-1}(s_i)B_{1,\chi\chi'}
\prod_{\ell\,|\,m}(1-\chi\chi'(\ell)). 
\end{align*}
If $\chi$ is nontrivial on $\langle p\rangle$, we have $\sum_{j=0}^{f-1}\chi(p^{-j})=0$, which implies that $\eqref{eq:tra2}=0$, i.e.,  \eqref{eq:tra1} holds. Hence, we assume that $\chi$
is trivial on $\langle p\rangle$. Then, the condition~\eqref{eq:tra1} is equivalent to  
\begin{equation}\label{eq:tra3}
\sum_{\chi'\in \AAA(E_0)}{\chi'}^{-1}(s_i)F_{\chi\chi'}=0, \quad 
\forall i=0,1,\ldots,e-1,\,  \forall \chi\in \D(m)\setminus \AAA(E_0) \mbox{ s.t. $\chi(p)=1$}, 
\end{equation}
where $F_{\chi\chi'}:=B_{1,\chi\chi'}
\prod_{\ell\,|\,m}(1-\chi\chi'(\ell))$. 

We now assume that \eqref{eq:tra3} holds. Then,  for any 
$\chi''\in \AAA(E_0)$, 
\begin{align*}
0=&\,\sum_{i=0}^{e-1}\sum_{\chi'\in \AAA(E_0)}{\chi'}^{-1}(s_i)F_{\chi\chi'}\chi''(s_i)\\
=&\,\sum_{\chi'\in \AAA(E_0)}F_{\chi\chi'}\left(\sum_{i=0}^{e-1}\chi''{\chi'}^{-1}(s_i)\right)
=eF_{\chi\chi''}.  
\end{align*}
Conversely, if $F_{\chi\chi''}=0$ for any  $\chi''\in \AAA(E_0)$, \eqref{eq:tra3} clearly holds. Hence, the condition~\eqref{eq:tra3} (i.e.,  the condition~\eqref{eq:tra1}) is equivalent to 
\begin{equation}\label{eq:tra4}
B_{1,\chi}
\prod_{\ell\,|\,m}(1-\chi(\ell))=0, \quad \forall \chi\in  \DD(m)\setminus \AAA(E_0) \mbox{ s.t. $\chi(p)=1$}. 
\end{equation}
If $\chi$ is an even character, $B_{1,\chi}=0$, which implies that 
\eqref{eq:tra4} holds. On the other hand, if  $\chi$ is an odd character, $B_{1,\chi}\not=0$, which implies that the condition~\eqref{eq:tra4} is equivalent to $\prod_{\ell\,|\,m}(1-\chi(\ell))=0$. This completes the proof. \qed
\vspace{0.3cm}

We now give a proof of Theorem~\ref{cor:equiv2}. 
\vspace{0.3cm}

{\bf Proof of Theorem~\ref{cor:equiv2}: }
Assume that $\mu_{p,f,m}=2$. By Lemma~\ref{prop:equiv1}, there is a subgroup $E_0$ of index $2$ of $G$ containing $H=\langle p\rangle$ such that for each $i=0,1$, $\sum_{j=0}^{f-1}[tp^j]_m$ is constant for any $t\in E_i$. 
Then, by Proposition~\ref{thm:mainpre},  
the condition (2) in the theorem holds. 
Since $A_0=A_1$ if and only if $\mu_{p,f,m}=1$ by Lemma~\ref{lem:pure1}, it must be $A_0\not=A_1$. We now see that $A_0\not=A_1$ if and only if the condition (1) holds. 
For $i=0,1$, we have    
\begin{align*}
A_i=&\,\frac{2f}{\phi(m)}\sum_{t\in E_i}t=\frac{f}{\phi(m)}\sum_{t\in G}t((-1)^i\chi_{\ann}(t)+1)
\\
=&\,\frac{(-1)^ifm}{\phi(m)}B_{1,\chi_{\ann}}\prod_{\ell\,|\,m}(1-\chi_\ann(\ell))+\frac{fm}{2}. 
\end{align*}
Hence, $A_0\not=A_1$ if and only if $B_{1,\chi_{\ann}}\prod_{\ell\,|\,m}(1-\chi_\ann(\ell))\not=0$, i.e., $\chi_{\ann}$ is an odd character and $\chi_\ann(\ell)=-1$ for any prime divisor $\ell$ of  $m$ not dividing the conductor of $\chi_{\ann}$ by noting that $\chi_{\ann}$ is of order $2$. 

Conversely, we assume that there is a subgroup $E_0$ of index $2$ of $G$ containing $H=\langle p\rangle$ satisfying the conditions (1) and (2). By Proposition~\ref{thm:mainpre}, for each $i=0,1$, $\sum_{j=0}^{f-1}[tp^j]_m$ is constant for any $t\in E_i$. Let $F$ be the  quadratic field fixed by the subgroup of $\Gal(K/\Q)$ corresponding to $E_0$. Now, we have $-1\not \in E_0$; otherwise, $\chi_{\ann}(-1)=1$, which contradicts to that $\chi_{\ann}$ is odd. This implies that $F$ is an imaginary quadratic field. 
On the other hand, by Theorem~\ref{thm:Stick}, we have 
\begin{align*}
(G_{p^f}(\eta_P^{-1})^m)O_{K}=P^{A_0\sum_{t\in E_0/H}\sigma_{t}+A_1\sum_{t\in E_0/H}\sigma_{s_1t}}.  
\end{align*}
Then, $P^{f\sum_{t\in E_0/H}\sigma_{t}}={P'}^d O_K$ for the ideal $P'=P\cap F$, where  $d=[O_K/P:O_F/P']$. Take a positive  integer $c$, 
e.g., the class number of $F$, so that 
${P'}^{cd}$ 
is principal in $F$. Then, ${P'}^{cd}=(\pi)$ for some integer $\pi$ in $F$. Since $-1\not \in E_0$, we have ${({P'}^{cd})}^{\sigma_{s_1}}=(\overline{\pi})$, where 
$\overline{\pi}$ is the complex conjugate of $\pi$. Furthermore, 
since $P^{f\sum_{t\in G/H}\sigma_{t}}=(p^f)$, we have 
$(\pi\overline{\pi})=(p^{cf})$. Thus, 
\[
(G_{p^f}(\eta_P^{-1})^{cfm})O_{K}=(\pi)^{A_0}(\overline{\pi})^{A_1}, 
\]
where $|\pi|=|\overline{\pi}|=p^\frac{cf}{2}$ and $A_0+A_1=fm$. 
Hence, $G_{p^f}(\eta_P^{-1})^{cfm}=\epsilon \pi^{A_0}\overline{\pi}^{A_1}$ for some unit $\epsilon$ of $K$ with $|\epsilon|=1$. Since $|\sigma(G_{p^f}(\eta_P^{-1})^{cfm})|=p^\frac{cf^2m}{2}$ for any $\sigma\in \Gal(K/\Q)$, it follows that 
$|\epsilon^\sigma|=1$. Hence, $\epsilon$ is a root of unity in $K$~(cf.~\cite[Theorem~2.1.13]{BEW97}). 
\qed
\vspace{0.3cm}

Note that the technique in the latter part of the proof above was also used in  previous papers, e.g., \cite{A10,FY,FYL,M98,M03,YLF,YX10}. 


\begin{remark}\label{rem:first1}
Assume that $G_{q}(\eta_P)^\ell\in F$ for some quadratic field $F$. 
Then, the following are obvious from Theorems~\ref{cor:equiv2} and \ref{thm:lift}. 
\begin{itemize}
\item[1.] $G_{q}(\eta_P^a)^\ell\in F$ for any $a$ with  $\gcd{(a,m)}=1$. 
\item[2.] For any $1\le i\le f-1$ with $\gcd{(i,f)}=1$, let $r$ be any prime such that $r\equiv p^i\,(\mod{m})$. Furthermore, let $\eta_m$ be a multiplicative character of order $m$ of $\F_{r^f}$. 
Then,  $G_{r^{f}}(\eta_m)^\ell$ is in some quadratic field.  
\item[3.] Let $s$ be any positive integer and $\eta$ be the lift of $\eta_P$ to $\F_{q^s}$. Then, $G_{q^{s}}(\eta)^\ell\in F$.   
\end{itemize}
\end{remark}

From Remark~\ref{rem:first1},  we define ${\mathcal P}$ (resp. ${\mathcal P}_f$) as the 
set of triples $(m,{\overline p},f)$ (resp. the set of pairs $(m,{\overline p})$) such that $\ord_{m}(p)=f$  and some nonzero integral power of $G_{p^f}(\eta_m)$ is in a quadratic field, where  ${\overline p}$ denotes a minimum representative of $\{[p^i]_m\mid 1\le i\le f-1,\gcd{(i,f)}=1\}$. Furthermore, we denote  the set of $(m,{\overline p},f)\in {\mathcal P}$ (resp. $(m,{\overline p})\in {\mathcal P}_f$) such that $A_0\not=A_1$ by ${\mathcal P}^\ast$ (resp. ${\mathcal P}_f^\ast$). 

Note that the Gauss sums in index $2$ case satisfy the condition in Proposition~\ref{prop:equiv1} as $E_0=\langle p\rangle$. 
Hence, defining ${\mathcal P}^{(2)}=\{(m,{\overline p},f) \mid \ord_{m}(p)=\phi(m)/2=f\}$, 
we have ${\mathcal P}^{(2)}\subseteq {\mathcal P}$. 


\section{Finiteness of ${\mathcal P}_{f}^\ast$}\label{sec:finite}
In this section, we prove the following theorem on the finiteness of  ${\mathcal P}_{f}^\ast$. 
\begin{theorem}\label{thm:finite}
The set ${\mathcal P}_{f}^\ast$ is finite for every positive integer $f$. 
\end{theorem}
To prove the theorem above, we prepare some notation. Let 
\begin{align*}
\OO^-(m,p):=&\,\{\chi\in \DD(m)\mid \chi(p)=1,\chi \mbox{ is an odd character}\}, \\
\XX^-(m,p):=&\,\{\chi\in \OO^-(m,p)\mid \mbox{The conductor of $\chi$ is divisible by any prime factor of $m$}\}. 
\end{align*}
Let $p_i$, $i=1,2,\ldots,r$, be distinct primes and $u_i$, $i=1,2,\ldots,r$, be
positive integers. 
Let $m=m_1m_2\cdots m_r$, where $m_i=p_i^{u_i}$, $i=1,2,\ldots,r$, and 
$f_i$, $i=1,2,\ldots,r$, denote the orders of $p$ modulo $m_i$, respectively. Then, $f=\lcm{(f_1,f_2,\ldots,f_r)}$. If $m$ is even, we will denote $m_1=2^{u_1}$ with $u_1\ge 1$. 

For odd $m_i$ or $m_i=m_1=2^{u_1}$ with $u_1=1,2$,  
let $\chi_i$ denote a fixed 
generator in $\DD({m_i})$. Note that $\chi_1$ is trivial if $m_1=2$.  
For $m_1=2^{u_1}$ with $u_1\ge 3$, since $(\Z/2^{u_1}\Z)^\times$ is generated by $-1$ and $5$, its character group is generated by $\chi_1'$ and $\chi_1''$ such that 
$\chi_1'(-1)=1$, $\chi_1'(5)=\zeta_{2^{u_1-2}}$, 
$\chi_1''(5)=1$ and $\chi_1''(-1)=-1$.  
Then, any character $\chi$ in $\DD(m)$ can be expressed as $\chi=\chi_1^{a_1}\chi_2^{a_2}\cdots \chi_r^{a_r}$ or 
$\chi={\chi_1'}^{a_0}{\chi_1''}^{a_1}\chi_2^{a_2}\cdots \chi_r^{a_r}$
for some integers $a_i$. 
For  odd $m_i$ or $m_i=m_1=4$, since $\chi_i^{f_i}(p)=\chi_i(p^{f_{i}})=1$ and $\chi_i^{f_i}(-1)=(-1)^{f_i}$, we have $\chi_{i}^{f_i}\in \OO^-(m,p)$ if and only if $f_i$ is odd. For $m_1=2^{u_1}$ with $u_1\ge 3$, let $f_1'$ and $f_1''$ be the orders of $p$ in $(\Z/m_1\Z)^\times/\langle -1\rangle$ and $(\Z/m_1\Z)^\times/\langle 5\rangle$. Then, ${\chi_1'}^{f_1'}{\chi_1''}^{f_1''}(p)=1$ and ${\chi_1'}^{f_1'}{\chi_1''}^{f_1''}(-1)=(-1)^{f_1''}$. 

By Theorem~\ref{cor:equiv2}, $\XX^-(m,p)=\{\chi_{\ann}\}$ or $\XX^-(m,p)=\emptyset$ if $(m,\overline{p},f)\in {\mathcal P}^\ast$. First, we treat the case where $\XX^-(m,p)=\{\chi_{\ann}\}$. 
\begin{proposition}\label{lem:finite1}
Assume that $(m,\overline{p})\in {\mathcal P}_{f}^\ast$ and $\XX^-(m,p)=\{\chi_{\ann}\}$. Then, $\phi(m_i)\,|\,4f$ holds for any $i=1,2,\ldots,r$. 
In particular, if $f$ is odd, it holds that $v_2(m)\in \{0,1,2,3\}$, where  $v_2(m)$ is the maximum integer such that $2^{v_2(m)}$ divides $m$. 
\end{proposition}
\pro 
For odd $m_i$ such that $f_i$ is even  and $f_i<\phi(m_i)/2$,  $\chi_i^{f_i}\chi_{\ann}\in \XX^{-}(m,p)$ and $\chi_i^{f_i}\chi_{\ann}\not=\chi_{\ann}$, a contradiction. Hence, $f_i\ge \phi(m_i)/2$. Since $f_i\,|\,\phi(m_i)$, it follows that  $f_i=\phi(m_i)$ or
$2f_i=\phi(m_i)$. Similarly, for odd $m_i$ such that  $f_i$ is odd and $f_i<\phi(m_i)/4$,  $\chi_i^{2f_i}\chi_{\ann}\in \XX^{-}(m,p)$ and  $\chi_i^{2f_i}\chi_{\ann}\not=\chi_{\ann}$, a contradiction. Hence, $2f_i=\phi(m_i)$ or
$4f_i=\phi(m_i)$  since $2f_i\,|\,\phi(m_i)$. Thus, in these cases, we have $\phi(m_i)\,|\,4f_i\,|\,4f$. 

For even $m_1=2^{u_1}$ with $u_1\ge 3$, we have $f_1'\ge 2^{u_1-3}$; otherwise, 
${\chi_1'}^{f_1'}\chi_{\ann}\in \XX^-(m,p)$ and ${\chi_1'}^{f_1'}\chi_{\ann}\not=\chi_{\ann}$, a contradiction. Since $f_1'$ is a power of $2$, we have  $m_1=2^{u_1}\,|\, 8f_1'\,|\,8f$, i.e., $\phi(m_1)\,|\,4f$ for any $u_1\ge 1$. 
\qed
\vspace{0.3cm}

The condition $\phi(m_i)\,|\,4f$ 
in Proposition~\ref{lem:finite1} implies that there are only  finite possibility for $m$ if $\XX^-(m,p)=\{\chi_{\ann}\}$.  Next, we treat the case where $\XX^-(m,p)=\emptyset$. 
\begin{proposition}\label{lem:finite2}
If $(m,\overline{p})\in {\mathcal P}_{f}^\ast$ and $\XX^-(m,p)=\emptyset$, $m$ is upper bounded by some positive integer determined by $f$. 
\end{proposition}
\pro 
%
For a subset $S$ of $\{1,2,\ldots,r\}$, let $m_S=\prod_{i\in S}m_i$. 
By Theorem~\ref{cor:equiv2}, there is $B\subseteq \{1,2,\ldots,r\}$ such that $m_B\,|\,m$ and $\chi_{\ann}\in \XX^-(m_B,p)$. Let $A$ 
be the set of $i\in \{1,2,\ldots,r\}$ such that $m_i$ is odd and  $4f_i\ge  \phi(m_i)$. We add $i=1$ into $A$ 
if $m_1\in \{2,4\}$ or $m_1=2^{u_1}$ with $u_1\ge 3$ and $4f_1'\ge  \phi(m_1)$. 
If $m_1=2^{u_1}$ with $u_1\in \{1,2,3\}$, it is clear that $1\in A$. 
Write $A_1=A\cap B$, $A_2=A\setminus B$, $B'=B\setminus A$ and $C=\{1,2,\ldots,r\}\setminus (A\cup B)$. 

Assume that $A_2=\emptyset$ and $C\not=\emptyset$. If $1\not\in C$ or 
$1\in C$ but $m_1$ is odd, we have $\chi_{\ann}\prod_{i\in C}\chi_i^{2f_i}\in \XX^-(m,p)$, a contradiction. If  $1\in C$ and $m_1=2^{u_1}$ with $u_1\ge 4$, we have $\chi_{\ann}{\chi_1'}^{f_1'}\prod_{i\in C\setminus \{1\}}\chi_i^{2f_i}\in \XX^-(m,p)$, a contradiction. Hence, if 
$A_2=\emptyset$, it follows that $C=\emptyset$, i.e., $\XX^-(m,p)=\{\chi_{\ann}\}$, a contradiction. Hence, $A_2\not=\emptyset$. 

Since $\phi(m_i)\le 4f_i\le 4f$ for any $i\in A$, there are only finite possibility
for $m_{A}$ and $A$ is finite.  Furthermore,  $1\le |A_2|\le |A|$. Let $C_1$ be the set of $i\in C$  
such that $m_i$ is odd and  $(2|A|+2)f_i\ge  \phi(m_i)$. We add $i=1\in C$ into $C_1$ if  
$m_1=2^{u_1}$ with $u_1\ge 4$ and $(2|A|+2)f_1'\ge  \phi(m_1)$.  
Furthermore, let $C_2=C\setminus C_1$. 
Note that there are only finite possibility for $m_{C_1}$.  

Next, we consider the  possibility for $m_{C_2}$.  
Firstly, assume that 
$m_1$ is odd or $m_1$ is even 
but $1\not\in C$. Then,  
by the pigeonhole principle, for any $k\in C_2$,  there are $j\in A_2$ and $h',h''\in \{1,2,\ldots,|A|+1\}$ with $h'<h''$ such that \[
(\chi_{\ann}\chi_k^{2h'f_k}\prod_{i\in C\setminus \{k\}}\chi_i^{2f_i})(p_j)=(\chi_{\ann}\chi_k^{2h''f_k}\prod_{i\in C\setminus \{k\}}\chi_i^{2f_i})(p_j)=1. 
\]
Hence, $\chi_k^{2(h''-h')f_k}(p_j)=1$, i.e., $p_j^{2(h''-h')f_k}\equiv 1\,(\mod{m_k})$. 
Secondly, assume that $m_1$ is even 
and $1(=k)\in C_2$. Then, it holds that  
\[
(\chi_{\ann}{\chi_1'}^{h'f_1'}\prod_{i\in C\setminus \{1\}}\chi_i^{2f_i})(p_j)=(\chi_{\ann}{\chi_1'}^{h''f_1'}\prod_{i\in C\setminus \{1\}}\chi_i^{2f_i})(p_j)=1. 
\]
Hence, ${\chi_1'}^{(h''-h')f_1'}(p_j)=1$, i.e., $p_j^{2(h''-h')f_1}\equiv 1\,(\mod{m_1})$. 
Thirdly, assume that $m_1$ is even 
and $1(\not=k)\in C$. Then,  
\[
(\chi_{\ann}\chi_k^{2h'f_k}{\chi_1'}^{f_1'}\prod_{i\in C\setminus \{k,1\}}\chi_i^{2f_i})(p_j)=(\chi_{\ann}\chi_k^{2h''f_k}{\chi_1'}^{f_1'}\prod_{i\in C\setminus \{k,1\}}\chi_i^{2f_i})(p_j)=1. 
\]
Hence, $\chi_k^{2(h''-h')f_k}(p_j)=1$, i.e., $p_j^{2(h''-h')f_k}\equiv 1\,(\mod{m_k})$. 
Thus, in all cases, $m_k$ for $k\in C_2$ is a divisor of $\prod_{j\in A_2}\prod_{h=1}^{|A|}(p_j^{2fh}-1)$. Hence, there are only finite possibility for   $m_{C_2}$. 

Finally, we consider the  possibility for $m_{B'}$. 
For any $i\in B'$ such that $m_i$ is odd, there is $j\in (A_2\cup C)=\{1,2,\ldots,r\}\setminus B$ such that $\chi_i^{2f_i}\chi_{\ann}(p_j)=1$. Since $\chi_{\ann}(p_j)=-1$, we have  $\chi_i^{2f_i}(p_j)=-1$, i.e., $p_j^{4f_i}\equiv 1\,(\mod{m_i})$. On the other hand, if $m_1$ is even 
and $1\in B'$, there is $j\in (A_2\cup C)=\{1,2,\ldots,r\}\setminus B$ such that ${\chi_1'}^{f_1'}\chi_{\ann}(p_j)=1$. Then, we have ${\chi_1'}^{f_1'}(p_j)=-1$, i.e., $p_j^{2f_1}\equiv 1\,(\mod{m_1})$.  
Thus, in all cases, $m_i$ for $i\in B'$  is a divisor of $\prod_{j\in (A_2\cup C)}(p_j^{4f}-1)$. Hence, there are only finite possibility for  $m_{B'}$. 
Thus, the assertion of the proposition follows. 
\qed

\vspace{0.3cm}
By Propositions~\ref{lem:finite1} and \ref{lem:finite2}, we obtain the assertion of  Theorem~\ref{thm:finite}. However, as in the proof of Proposition~\ref{lem:finite2}, $m$ is upper bounded by a very large number. On the other hand, if $f$ is odd, we have a sharper estimation for $m$. 

\begin{theorem}\label{thm:fodd}
Assume that $f$ is odd,  $(m,\overline{p})\in {\mathcal P}_{f}^\ast$ and $\XX^-(m,p)=\emptyset$.  Then,   $v_2(m)\le 2$ and the following  hold. 
\begin{itemize}
\item[(1)] If $v_2(m)\in \{0,2\}$, $r$ is even and $f_i=\phi(m_i)/2$ for every $i=1,2,3\ldots,r$. 
\item[(2)] If $v_2(m)=1$ and $r$ is even,  we have either  $f_j=\phi(m_j)/2$, $m_j\,|\, 2^{2f_j}-1$, or $m_j\,|\, 2^{4f_j}-1$ and $6f_j=\phi(m_j)$ for any $i=2,3,\ldots,r$. 
\item[(3)] If $v_2(m)=1$ and $r$ is odd, we have  either  $f_j=\phi(m_j)/4$, $f_j=\phi(m_j)/2$, $m_j\,|\, 2^{2f_j}-1$, or $m_j\,|\, 2^{4f_j}-1$ and $8f_j=\phi(m_j)$  for any $i=2,3,\ldots,r$. 
\end{itemize}
\end{theorem}
\pro
We first assume that $m$ is odd. 
Define $
\theta=\prod_{i=1}^r\chi_i^{f_i}$. Since $f$ is odd,  $r$ is  even; otherwise $\theta(-1)=-1$, which is a contradiction to $\XX^-(m,p)=\emptyset$. 
Furthermore, if $f_i<\phi(m_i)/2$ for some $i$, we have  $\chi_i^{f_i}\theta\in \XX^{-}(m,p)$, a contradiction. Hence, $\phi(m_i)=2f_i\,|\,2f$ for every $i$. 

We next assume that $m$ is even. Let $m_1=2^{u_1}$ with $u_1\ge 1$. If $u_1\ge 3$, $f_1=1$ and either of $\chi_1'\prod_{i=2}^r\chi_i^{f_i}$ or $\chi_1''\prod_{i=2}^r\chi_i^{f_i}$ is in $\XX^{-}(m,p)$, a contradiction. Hence, we have $u_1=1$ or $2$.  
If $u_1=2$, similarly to the case where $m$ is odd, $r$ must be even and $\phi(m_i)=2f_i\,|\,2f$ for every $i$. Hence, the assertion (1) of the theorem follows. 

Finally, we assume that $u_1=1$, i.e.,  $m=2m_2m_3\cdots m_r$. Define 
$\theta'=\prod_{i=2}^r\chi_i^{f_i}$. 
If $r$ is even, $\theta'$ and $\chi_j^{2f_j}\theta'$ are in $\XX^-(m/2,p)$ for any $j$ such that $f_j<\phi(m_j)/3$. Hence, by Theorem~\ref{cor:equiv2}, we have 
$\theta'(2),\chi_j^{2f_j}\theta'(2)\in \{1,-1\}$, i.e., $\chi_j^{2f_j}(2)=1$ or $-1$.  
If $\chi_j^{2f_j}(2)=1$, we have $m_j\,|\, 2^{2f_j}-1$. If   $\chi_j^{2f_j}(2)=-1$,  we have 
 $\chi_j^{2f_j}\theta'=\chi_{\ann}$. Hence, $m_j\,|\, 2^{4f_j}-1$ and $6f_j=\phi(m_j)$ hold. Therefore, the assertion (2) of  the theorem follows. 
If $r$ is odd, $\chi_j^{f_j}\theta'$ and $\chi_j^{3f_j}\theta'$ are in $\XX^-(m/2,p)$ for any $j$ such that $f_j<\phi(m_j)/4$. Hence, by Theorem~\ref{cor:equiv2}, we have 
$\chi_j^{f_j}\theta'(2),\chi_j^{3f_j}\theta'(2)\in \{1,-1\}$, i.e., $\chi_j^{2f_j}(2)=1$ or $-1$. If $\chi_j^{2f_j}(2)=1$, we have $m_j\,|\, 2^{2f_j}-1$. If   $\chi_j^{2f_j}(2)=-1$,  we have 
 $\chi_j^{3f_j}\theta'=\chi_{\ann}$. Hence, $m_j\,|\, 2^{4f_j}-1$ and $8f_j=\phi(m_j)$ hold.  Therefore, the assertion (3) of the theorem follows. 
\qed
\vspace{0.3cm}

By Proposition~\ref{lem:finite1} and Theorem~\ref{thm:fodd}, we obtain the following corollary, which corresponds to \cite[Corollary~5.2]{A04} for pure Gauss sums. 
\begin{corollary}\label{cor:fodd2}
Assume that  $f$ is odd and $(m,\overline{p})\in {\mathcal P}_f^\ast$. Then, $v_2(m)\in \{0,1,2,3\}$ and $\frac{m}{2^{v_2(m)}}\,|\,2^{4f}-1$.   
\end{corollary}
\pro
Both $f_j=\phi(m_j)/2$ and $f_j=\phi(m_j)/4$ imply that  
$m_j\,|\, 2^{4f_j}-1$ by Euler's theorem. 
\qed
\vspace{0.3cm}

Now, we can classify  $(m,\overline{p})\in {\mathcal P}_f^\ast$ for small 
odd $f$. 
\begin{corollary}
The sets ${\mathcal P}_f^\ast$ for $f=1,3,5,7$ are determined as follows: 
\begin{align*}
{\mathcal P}_1^\ast&\,=\{(3,1),(4,1),(6,1)\};\\
{\mathcal P}_3^\ast&\,=\{(7,2),(9,4),(18,7),(21,4),(28,9),(39,16)\};\\
{\mathcal P}_5^\ast&\,=\{(11,3),(22,3),(33,4),(55,16),(66,25)\};\\
{\mathcal P}_7^\ast&\,=\emptyset. 
\end{align*} 
\end{corollary}
\pro 
First, we list all $m$ satisfying the condition of Corollary~\ref{cor:fodd2}. 
Then,  we reduce the  candidates of $m$ such that $(m,\overline{p})\in {\mathcal P}_f^\ast$ by applying Proposition~\ref{lem:finite1} and Theorem~\ref{thm:fodd}.  
For remaining candidates, we used a computer to  directly check whether there is $\overline{p}$ such that $(m,\overline{p})\in {\mathcal P}_f^\ast$ based on Lemma~\ref{prop:equiv1} for $\mu_{p,f,m}=2$.  
\qed

\section{Characterization of $(m,\overline{p},f)\in {\mathcal P}^\ast$ in the $r\le 2$ case}\label{sec:m1m2}
In this section,  we assume that $m$ has at most two distinct prime power divisors, i.e., $m=p_1^{u_1}$ or $m=p_1^{u_1}p_2^{u_2}$. In particular, when 
$m=p_1^{u_1}p_2^{u_2}$, we assume that each $m_i$ is odd or $m_1=2^{u_1}$ with $u_1\in \{1,2\}$. The objective of this section is to 
 characterize $(m,\overline{p},f)\in {\mathcal P}^\ast$ in these cases. Note that 
any character in $G$ has the form $\chi_1^{a_1}$,  ${\chi_1'}^{a_0}{\chi_1''}^{a_1}$ or $\chi_1^{a_1}\chi_2^{a_2}$ in these cases. 
Though we are also able to characterize $(m=2^{u_1}m_2,\overline{p},f)\in {\mathcal P}^\ast$ for $u_1\ge 3$, the argument becomes more complicated since each character in $G$ has the form ${\chi_1'}^{a_0}{\chi_1''}^{a_1}\chi_2^{a_2}$.  Hence, we will not treat this case. 
\subsection{The $r=1$ case}
In this subsection, we characterize $(m=p_1^{u_1},\overline{p},f)\in {\mathcal P}^\ast$. 
\begin{theorem}\label{prop:fac1}
Assume that $m=m_1=p_1^{u_1}$ with $p_1$ odd or $p_1=2$ and $u_1\in \{1,2\}$. Then, $(m,\overline{p},f)\in {\mathcal P}^\ast$ if and only if  $(m,\overline{p},f)=(4,1,1)$ or  $p_1\equiv 3\,(\mod{4})$ and $f=f_1=\phi(m_1)/2$.  
\end{theorem}
\pro
%
%
Since $\chi_{\ann}$ is of order $2$ and an odd character modulo $m_1$, 
$\phi(m_1)/2$ is odd, i.e., $p_1\equiv 3\,(\mod{4})$ or $m_1=4$.  
If $m_1=4$, we have $f=f_1=1$ and $\overline{p}=1$. 
We now assume that $p_1\equiv 3\,(\mod{4})$. Since $\chi_{\ann}(p)=1$,  $f_1$ is a divisor of 
$\phi(m_1)/2$, This implies that $f_1$ is odd.  
If $f=f_1< \phi(m_1)/4$, $\chi_1^{2f_1}\chi_{\ann}\in \XX^-(m_1,p)$ and  $\chi_1^{2f_1}\chi_{\ann}\not=\chi_{\ann}$, a contradiction. Hence, either  $f_1=\phi(m_1)/2$  or   $f_1=\phi(m_1)/4$ holds. But,  $f_1=\phi(m_1)/4$ is 
impossible since $\phi(m_1)/2$ is odd.  Thus, we obtain $f_1=\phi(m_1)/2$. 

The converse is true  due to Theorem~\ref{cor:equiv2} since  $\D^-(m,p)=\{\chi_1^{\phi(m_1)/2}\}$. 
\qed

\begin{theorem}\label{prop:fac2}
Assume that $m=m_1=p_1^{u_1}$ with $p_1=2$ and $u_1\ge 3$. Then, $(m,\overline{p},f)\in {\mathcal P}^\ast$ if and only if  $f_1'=2^{u_1-2}$. 
\end{theorem}
\pro
%
%
Since $\chi_{\ann}$ is of order $2$ and an odd character, we have $\chi_{\ann}=\chi_1''$ or $\chi_{\ann}={\chi_1'}^{2^{u_1-3}}\chi_1''$. If 
$f_1'<2^{u_1-2}$, we have ${\chi_1'}^{f_1'}\chi_{\ann}\in \XX^-(m,p)$ and ${\chi_1'}^{f_1'}\chi_{\ann}\not=\chi_{\ann}$,  a contradiction. Hence, $f_1'=2^{u_1-2}$. 

The converse is true  due to Theorem~\ref{cor:equiv2} since  $\D^-(m,p)=\{\chi_1''\}$ or $\{{\chi_1'}^{2^{u_1-3}}\chi_1''\}$ depending on whether $f_1''=1$ or $2$, respectively. 
\qed
\vspace{0.3cm}

From Theorems~\ref{prop:fac1} and \ref{prop:fac2}, all $(m=m_1,\overline{p},f)\in {\mathcal P}^\ast$ are in the index $2$ case. 
\subsection{The  $r=2$ case}
In this subsection, we characterize $(m=m_1m_2,\overline{p},f)\in {\mathcal P}^\ast$.  

For two positive integers $d_1,d_2$, 
let 
\[
A^-(d_1,d_2)=\left\{(a_1,a_2)\in \Z^2\mid \frac{a_1}{d_1}+\frac{a_2}{d_2}\in \Z,0<a_i<d_i\, \,  (i=1,2), a_1+a_2\equiv 1\,(\mod{2})\right\}. 
\] 
Similarly, we define $A^+(d_1,d_2)$ by replacing $a_1+a_2\equiv 1\,(\mod{2})$ in the definition of $A^-(d_1,d_2)$ by  $a_1+a_2\equiv 0\,(\mod{2})$. Observe that when either one of $d_i$, $i=1,2$, is odd, 
$A^-(d_1,d_2)=\emptyset$ if and only if $A^+(d_1,d_2)=\emptyset$ since  $(a_1,a_2)\in A^-(d_1,d_2)$ if and only if $(d_1-a_1,d_2-a_2)\in A^+(d_1,d_2)$. Furthermore, when $d_1$ and $d_2$ are both odd,  $A^+(d_1,d_2)=\emptyset$. 

We can choose $\chi_i$ so that $\chi_i(p)=\zeta_f^{f/f_i}$. 
For $\chi\in \D(m)$, we write $\chi=\chi_1^{a_1}\chi_2^{a_2}$  for some integers $a_i$. Then, we have 
\[
\chi(p)=\zeta_f^{f(\frac{a_1}{f_1}+\frac{a_2}{f_2})} \mbox{ and } 
\chi(-1)=(-1)^{a_1+a_2}. 
\] 
Hence, we have the following lemma.  
\begin{lemma}\label{lem:aokidd2}
Let $m=m_1m_2$ with $m_i$ odd or $m_1=2^{u_1}$ with $u_1\in \{1,2\}$. 
Assume that  $(m,\overline{p},f)\in {\mathcal P}^\ast$. Then, $\XX^{-}(m,p)$ consists of only characters of the form $\chi_1^{i_1f_1}\chi_2^{i_2f_2}$ for some integers $i_1,i_2$ if and only if $A^-(f_1,f_2)=\emptyset$.  
\end{lemma}
%
In \cite{A04}, the following is known. 
\begin{lemma}{\em (\cite[Lemma 7.1]{A04})} \label{lem:aokidd}
The following hold. 
\begin{itemize}
\item[(i)] Assume  that $\gcd{(d_1,d_2)}$ is even. Then, $A^-(d_1,d_2)=\emptyset$ if and only if $v_2(d_1)=v_2(d_2)$. 
\item[(ii)] Assume  that $\gcd{(d_1,d_2)}$ is odd. Then, $A^-(d_1,d_2)=\emptyset$ if and only if $\gcd{(d_1,d_2)}=1$. 
\end{itemize}
\end{lemma}
By using Lemma~\ref{lem:aokidd}, we have the following.  
\begin{lemma}\label{lem:f1f2half}
Assume that $f_1,f_2\equiv 0\,(\mod{2})$, $f_1/2\equiv 1\,(\mod{2})$ and  
$f_2/2\equiv 0\,(\mod{2})$. Let $g_1=f_1/2$ and $g_2=f_2/2$.  Then, 
$A^-(f_1,f_2)=\{(g_1,g_2)\}$ if and only if $\gcd{(g_1,g_2)}=1$. 
\end{lemma} 
\pro
Define 
\begin{align*}
S_1=&\,\left\{(a_1,a_2)\in \Z^2\mid \frac{a_1}{g_1}+\frac{a_2}{g_2}\in 2\Z,0< a_i<f_i, a_i\not=g_i\, (i=1,2), a_1+a_2\equiv 1\,(\mod{2})\right\},\\  
S_2=&\,\left\{(a_1,a_2)\in \Z^2\mid \frac{a_1}{g_1}+\frac{a_2}{g_2}\in \Z,0<a_i<f_i, a_i\not=g_i\, (i=1,2), a_1+a_2\equiv 1\,(\mod{2})\right\},\\
S_3=&\, \left\{(a_1,a_2)\in \Z^2\mid\frac{a_1}{g_1}+\frac{a_2}{g_2}\in \Z,\, 0< a_i<g_i\, (i=1,2)\right\}. 
\end{align*}

Assume that $A^-(f_1,f_2)=\{(g_1,g_2)\}$. This condition is equivalent to  
$S_1=\emptyset$. 
Furthermore, this implies that  $S_2=\emptyset$. 
In fact, if there is $(a_1,a_2)\in S_2$,  
either $(a_1,a_2)$, $(a_1,a_2-g_2)$ or $(a_1,a_2+g_2)$ is in $S_1$. Then, by restricting the range $0< a_i<f_i$ of $S_2$ to  $0< a_i<g_i$ for $i=1,2$, we obtain $A^-(g_1,g_2)=\emptyset$. Finally, 
by Lemma~\ref{lem:aokidd}, it holds that  $\gcd{(g_1,g_2)}=1$. 

Conversely, assume that $\gcd{(g_1,g_2)}=1$. Then, by the Chinese remainder theorem, $S_3=\emptyset$.  
Furthermore, this implies that $S_1=\emptyset$. 
In fact, if there is  $(a_1,a_2)\in S_1$, 
either $(a_1,a_2)$, $(a_1,a_2-g_2)$, $(a_1-g_1,a_2)$ or $(a_1-g_1,a_2-g_2)$ is in  $S_3$.  Hence,  $A^-(f_1,f_2)=\{(g_1,g_2)\}$ follows. 
\qed
\vspace{0.3cm}

We now  characterize $(m=m_1m_2,\overline{p},f)\in {\mathcal P}^\ast$. We first consider the case where $\XX^-(m,p)=\{\chi_{\ann}\}$. In the following theorem, we  allow that $m_1$ is odd and $m_2$ is even for convenience. 
\begin{theorem}\label{prop:m1m2}
Assume that $m=m_1m_2$ with  $v_2(m)\in \{0,1,2\}$. Then, $(m,\overline{p},f)\in {\mathcal P}^\ast$ and $\XX^-(m,p)=\{\chi_{\ann}\}$ if and only if either of the following holds. 
%
\begin{itemize}
\item[(1)] $f_1=\phi(m_1)/2$ and $f_2=\phi(m_2)/4$ are odd, $\gcd{(f_1,f_2)}=1$, $p_2$ is quadratic modulo $p_1$, and $p_1$ is quartic  modulo $p_2$. 
\item[(2)] $f_1=\phi(m_1)/2$ is odd, $f_2=\phi(m_2)/2$ is even, $\gcd{(f_1,f_2)}=1$ and $p_2$ is quadratic modulo $p_1$. 
\item[(3)] $f_1=\phi(m_1)$, $f_2=\phi(m_2)$, $f_1/2$ is odd, $f_2/2$ is even  and $\gcd{(f_1/2,f_2/2)}=1$.  
\end{itemize}
\end{theorem}
\pro 
First, assume that $(m,\overline{p},f)\in {\mathcal P}^\ast$ and $\XX^-(m,p)=\{\chi_{\ann}\}$. Consider the following three  cases: 
(i) both $f_1$ and $f_2$ are odd; (ii)  $f_1$ is odd and $f_2$ is even; 
(iii) both $f_1$  and $f_2$ are even. Note that $v_2(m)=1$ is impossible since 
$X^-(m,p)=\{\chi_{\ann}\}$. 

{\bf Case (i):} If $f_i<\phi(m_i)/4$ for some $i=1,2$, we have $\chi_i^{2f_i}\chi_{\ann}\in \XX^-(m_1m_2,p)$ and $\chi_i^{2f_i}\chi_{\ann}\not=\chi_{\ann}$, a contradiction. Hence, 
$f_i\ge \phi(m_i)/4$. In particular, since $2f_i\,|\,\phi(m_i)$, 
either $f_i= \phi(m_i)/2$ or $f_i= \phi(m_i)/4$ follows. On the other hand, since 
$\chi_{\ann}=\chi_1^{\phi(m_1)/2}\chi_2^{\phi(m_2)/2}\in \XX^-(m_1m_2,p)$,  one of 
$\phi(m_i)/2$ is odd and another is even. Hence, we can assume that  
$f_1= \phi(m_1)/2$ and $f_2= \phi(m_2)/4$. Furthermore, by Lemmas~\ref{lem:aokidd2} and \ref{lem:aokidd},  $\gcd{(f_1,f_2)}=1$ follows. Finally, we have 
$\chi_1^{f_1}(p_2)=\chi_2^{f_2}(p_1)=1$ by Theorem~\ref{cor:equiv2}. 
This shows the assertion (1). 

{\bf Case (ii):} If $f_1<\phi(m_1)/4$, we have $\chi_1^{2f_1}\chi_{\ann}\in \XX^-(m_1m_2,p)$, a contradiction. Hence, 
$f_1= \phi(m_1)/2$ or $f_1=\phi(m_1)/4$. Similarly,   if $f_2<\phi(m_2)/2$, we have $\chi_2^{f_2}\chi_{\ann}\in \XX^-(m_1m_2,p)$, a contradiction. Hence, 
$f_2= \phi(m_2)$ or $f_2=\phi(m_2)/2$.  
If $f_2=\phi(m_2)$, since $\chi_{\ann}(p)=\chi_1^{\phi(m_1)/2}\chi_2^{\phi(m_2)/2}(p)=1$, we have $\chi_2^{\phi(m_2)/2}(p)=\chi_1^{\phi(m_1)/2}(p)=-1$. Hence, both of $f_1= \phi(m_1)/2$ and  $f_1=\phi(m_1)/4$ are impossible.  If $f_2=\phi(m_2)/2$, since 
$\chi_{\ann}(-1)=\chi_1^{\phi(m_1)/2}\chi_2^{\phi(m_2)/2}(-1)=-1$, we have $\chi_1^{\phi(m_1)/2}(-1)=-1$, which implies that 
$f_1=\phi(m_1)/2$. Furthermore, by Lemmas~\ref{lem:aokidd2} and \ref{lem:aokidd}, it follows that $\gcd{(f_1,f_2)}=1$. Finally,  we have  $\chi_1^{f_1}(p_2)=1$ by Theorem~\ref{cor:equiv2}. This shows the assertion (2). 

{\bf Case (iii):} If $f_i<\phi(m_i)/2$ for some $i=1,2$, we have $\chi_i^{f_i}\chi_{\ann}\in \XX^-(m_1m_2,p)$ and $\chi_i^{f_i}\chi_{\ann}\not=\chi_{\ann}$, a contradiction. Hence, either 
$f_i=\phi(m_i)/2$ or $f_i=\phi(m_i)$ follows. 
On the other hand, since $\chi_{\ann}(-1)=\chi_1^{\phi(m_1)/2}\chi_2^{\phi(m_2)/2}(-1)=-1$, one of $\phi(m_1)/2$ or $\phi(m_2)/2$ is odd and another is even. We can assume that $\phi(m_1)/2$ is odd and  $\phi(m_2)/2$ is even. Hence, we have $f_1=\phi(m_1)$. Furthermore, since $\chi_1^{\phi(m_1)/2}(p)=-1$, we have $\chi_2^{\phi(m_2)/2}(p)=-1$, i.e., 
$f_2=\phi(m_2)$. Since 
$A^-(f_1,f_2)=\{(f_1/2,f_2/2)\}$,  it follows that
$\gcd{(f_1/2,f_2/2)}=1$ from Lemma~\ref{lem:f1f2half}. This shows the assertion (3). 

Finally, the converses are true due to Theorem~\ref{cor:equiv2}  since $\D^-(m,p)=\{\chi_1^{f_1}\chi_2^{2f_2},\chi_1^{f_1},\chi_2^{f_2}\}$, $\D^-(m,p)=\{\chi_1^{f_1}\chi_2^{f_2},\chi_1^{f_1}\}$ and $\D^-(m,p)=\{\chi_1^{f_1/2}\chi_2^{f_2/2}\}$ in the cases (1), (2) and (3), respectively.  \qed

\begin{remark}\label{rem:notriv}
The Gauss sums in the case (1) of Theorem~\ref{prop:m1m2} are out of the framework of previous studies in literature. In fact, $\phi(m)/f=8$ and the group $(\Z/m\Z)^\times/\langle p\rangle $ is not elementary abelian. Hence, 
we have infinitely many examples of Gauss sums, some nonzero integral powers of which are in quadratic fields but not belong to the index $2$ case. For example, let $m_1=p_1\equiv 1\,(\mod{10})$ and $m_2=5$. It is clear that 
$\gcd{(\phi(m_1)/2,\phi(m_2)/4)}=1$ and $p_1$ is quartic modulo $p_2$. Furthermore, by the quadratic reciprocity law, $p_2$ is quadratic modulo $p_1$. Take a prime $p$ such  that $p\equiv 1\,(\mod{5})$ and $p\equiv g^2\,(\mod{p_1})$ for a primitive element $g$ of $\F_{p_1}$. Then, $(m,\overline{p},f)\in {\mathcal P}^\ast$, which is in the case (1) of Theorem~\ref{prop:m1m2}.  
\end{remark}
We next consider the case where $\XX^-(m,p)=\emptyset$. 
\begin{theorem}
Assume that $m=2m_2=2p_2^{u_2}$ with $p_2$ odd. 
Then, $(m,\overline{p},f)\in {\mathcal P}^\ast$ if and only if $f=f_2=\phi(m_2)/2$ is odd and $p_2\equiv 3\,(\mod{8})$.  
\end{theorem}
\pro 
Assume that $(m,\overline{p},f)\in {\mathcal P}^\ast$.  It is clear that 
$\XX^-(m,p)=\emptyset$. 
Since $\chi_{\ann}=\chi_{2}^{\phi(m_2)/2}$, $\phi(m_2)/2$ is odd and $f=f_2\,|\,\phi(m_2)/2$. Since $\chi_{\ann}(2)=-1$, it follows that  $p_2\equiv 3\,(\mod{8})$ by the supplementary law of quadratic reciprocity. Furthermore, $f_2$ is odd. If $f_2<\phi(m_2)/2$, we have $\chi_2^{f_2}(2)=1$ but  $\chi_2^{\phi(m_2)/2}(2)=-1$, a contradiction. Therefore, $f_2=\phi(m_2)/2$. 

The converse also holds due to Theorem~\ref{cor:equiv2} since 
$\OO^-(m,p)=\{\chi_2^{f_2}\}$. 
\qed
\vspace{0.3cm}

Similarly to Theorem~\ref{prop:m1m2}, we  allow that  $m_1$ is odd and $m_2$ is even in the following theorem. 
\begin{theorem}
Assume that $m=m_1
m_2=p_1^{u_1}p_2^{u_2}$ with $v_2(m)\in \{0,2\}$. Then, $(m,\overline{p},f)\in {\mathcal P}^\ast$ and $\XX^-(m,p)=\emptyset$ if and only if either of the following holds. 
\begin{itemize}
\item[(1)] $f_1=\phi(m_1)/2$ and $f_2=\phi(m_2)/2$ are odd, $\gcd{(f_1,f_2)}=1$,  $p_2$ is nonquadratic modulo $p_1$,  and
$p_1$ is quadratic modulo $p_2$. 
\item[(2)] $f_1=\phi(m_1)/2$ is odd, $f_2=\phi(m_2)$, $\gcd{(f_1,f_2)}=1$ and $p_2$ is nonquadratic modulo $p_1$.  
\end{itemize}
\end{theorem}
\pro 
First, assume that $(m,\overline{p},f)\in {\mathcal P}^\ast$ and $\XX^-(m,p)=\emptyset$.  Consider the following three  cases: 
(i) both $f_1$ and $f_2$ are odd; (ii)  $f_1$ is odd and $f_2$ is even; 
(iii) both $f_1$  and $f_2$ are even. 

{\bf Case (i):} 
If $f_1<\phi(m_1)/2$, we have $\chi_1^{2f_1}\chi_2^{f_2}\in \XX^{-}(m,p)$, a contradiction. Hence, $f_1=\phi(m_1)/2$. Similarly,  $f_2=\phi(m_2)/2$. Here, we can assume that  $\chi_{\ann}=\chi_1^{f_1}$. Then,  we have  
$\chi_1^{f_1}(p_2)=-1$ and $\chi_2^{f_2}(p_1)=1$ since $\chi_i^{f_i}\in \D^-(m,p)$.  
Furthermore, by Lemmas~\ref{lem:aokidd2} and \ref{lem:aokidd}, it follows that $\gcd{(f_1,f_2)}=1$. This shows the assertion (1). 

{\bf Case (ii):} If $f_2\not=\phi(m_2)$, we have $\chi_1^{f_1}\chi_2^{f_2}\in \XX^{-}(m,p)$, a contradiction. Hence, $f_2=\phi(m_2)$. Then, the characters in $\OO^-(m,p)$ are only $\chi_{1}^{if_1}$ for odd $i$. Hence, $\chi_{\ann}=\chi_{1}^{jf_1}$ for some odd $j$. If $\chi_1^{f_1}(p_2)=1$, we have $\chi_{\ann}(p_2)=1$, a contradiction. This implies that $\chi_{\ann}=\chi_1^{f_1}$ and $\chi_1^{f_1}(p_2)=-1$. Since $\chi_{\ann}$ is of order $2$, we have 
$f_1=\phi(m_1)/2$. Furthermore, by Lemmas~\ref{lem:aokidd2} and \ref{lem:aokidd}, it follows that $\gcd{(f_1,f_2)}=1$. This shows the assertion (2).

{\bf Case (iii):} Since $f_1$ and $f_2$ are even, we have $p^{f_i/2}\equiv -1\,(\mod{m_i})$ for $i=1,2$.  Furthermore, by Lemmas~\ref{lem:aokidd2} and \ref{lem:aokidd}, it follows that $v_2(f_1)=v_2(f_2)=v_2(f)$. Then, we obtain 
 $p^{f/2}\equiv -1\,(\mod{m})$, i.e., $(m,\overline{p},f)$ is in the semi-primitive case.  
Hence, this case is impossible. 

Finally, the converses are true due to Theorem~\ref{cor:equiv2} since $\OO^-(m,p)=\{\chi_1^{f_1},\chi_2^{f_2}\}$ and $\OO^-(m,p)=\{\chi_1^{f_1}\}$ in the cases (1) and (2), respectively.  
\qed

 \section{Concluding remarks}
In this paper, we gave a general necessary condition for some nonzero integral power of a  Gauss sum to be in an extension field of  minimal degree $e$ over $\Q$ for a fixed $e\ge 1$ in terms of Dirichlet characters modulo $m$. In particular, we proved that the set of pairs $(m,\overline{p})$ such that $\ord_{m}(p)=f$ and some nonzero integral power of $ G_{p^f}(\eta_m)$ is in a quadratic field $F/\Q$ is finite for each $f$. 
However, we could not prove this claim for extension fields of  degree  $e>2$. Hence, 
the following problem naturally arises. 
\begin{problem}
Determine whether the set of pairs $(m,\overline{p})$ such that $\ord_{m}(p)=f$ and some nonzero integral power of $ G_{p^f}(\eta_m)$ is in  an   extension field of minimal degree $e$ over $\Q$ is finite for fixed $e>2$ and $f\ge 1$. 
\end{problem}
We give a list of $(m,{\overline p},f)\in {\mathcal P}^\ast$ for $m\le 1000$ except for those in index $2$ case in Tables~\ref{Tab0} and \ref{Tab1}. The computation is based on the claim of Lemma~\ref{prop:equiv1} for $\mu_{p,f,m}=2$ and its converse shown in the proof of Theorem~\ref{cor:equiv2} as follows. For each 
$m\le 1000$ and $1\le \overline{p}\le m-1$ such that $\gcd{(\overline{p},m)}=1$, check whether $\sum_{i=0}^{f-1}[x\overline{p}^i]_m$  for representatives $x$ of $(\Z/m\Z)^\times/\langle \overline{p}\rangle $ take exactly two values, where $f$ is the order of $\overline{p}$ in  $(\Z/m\Z)^\times$. If the set of $x\in (\Z/m\Z)^\times$ such that $\sum_{i=0}^{f-1}[xp^i]_m=\sum_{i=0}^{f-1}[p^i]_m$ forms a subgroup of index $2$ in $(\Z/m\Z)^\times$, then return $[m,\overline{p},f,|(\Z/m\Z)^\times/\langle \overline{p}\rangle |]$ except for the case  where $\langle \overline{p}\rangle$ is of index $2$ in $(\Z/m\Z)^\times$. 
As seen from the tables, the index $|(\Z/m\Z)^\times/\langle p\rangle|$ is a power of $2$ for all listed examples. Indeed, this is correct  when $f$ is odd and $v_2(m)\not=1$ by Proposition~\ref{lem:finite1} and Theorem~\ref{thm:fodd} or when $m$ has at most two distinct prime divisors by the results in Section~\ref{sec:m1m2}. 
However,  it remains unsolved in general. Furthermore,  one may be interested in the structure of $(\Z/m\Z)^\times/\langle p\rangle$. Hence, we give the following open problem.
\begin{problem}
Determine the structure of the group $(\Z/m\Z)^\times/\langle p\rangle$ for  $(m,\overline{p},f)\in {\mathcal P}^\ast$. 
\end{problem}
When $(\Z/m\Z)^\times/\langle p\rangle$ is an elementary abelian $2$-group, the corresponding Gauss sums were already evaluated in \cite{A10}. 
However,  $(\Z/m\Z)^\times/\langle p\rangle$ for $(m,\overline{p},f)\in {\mathcal P}^\ast$ is not necessarily elementary abelian as seen in Remark~\ref{rem:notriv}. Thus, the problem treated in this paper is out of the framework of the study in \cite{A10}. 

The condition for $(m,p)$ such that $\XX^-(m,p)=\emptyset$ was studied in \cite{A04,A12} in relation to Diophantine equations 
of the form  $\sum_{i=1}^nx_i/d_i\equiv 0\,(\mod{1})$. On the other hand, in Proposition~\ref{lem:finite2}, 
we gave a short proof for that the set of $(m,\overline{p})\in {\mathcal P}_{f}^\ast$ such that $\XX^-(m,p)=\emptyset$ is finite without using this argument. 
One may use the results in \cite{A04,A12} to derive a good bound 
for $m$ in the case where $f$ is even.  
Hence, we give the following problem for future work.
\begin{problem}
Give a sharp upper bound for $m$ such that $(m,\overline{p})\in {\mathcal P}_f^\ast$ when $f$ is even. Furthermore, classify $(m,\overline{p})\in {\mathcal P}_f^\ast$ for small even $f$. 
\end{problem}

Finally, we give one more important problem. 
We studied the triples $(m,\overline{p},f)\in {\mathcal P}^\ast$ but did not give an explicit evaluation for the corresponding Gauss sums. To do this, we need to determine the minimum $\ell$ such that  $G_{p^f}(\eta_m)^\ell\in F$. Hence, we give the following problem. 
\begin{problem}\label{prob:last}
For $(m,\overline{p},f)\in {\mathcal P}^\ast$, determine the minimum $\ell $ such that  $G_{p^f}(\eta_m)^\ell$ is in some quadratic field $F$. Furthermore, explicitly evaluate $G_{p^f}(\eta_m)$ and $G_{p^f}(\eta_m)^\ell$ in this case.  
\end{problem}


\newpage
{\tiny 
\begin{table}[h]
\caption{The list of $(m,{\overline p},f)\in {\mathcal P}^\ast$ for $m\le 1000$ under the assumption that  $(h:=)|(\Z/m\Z)^\times/\langle p\rangle|> 2$ }
\label{Tab0}
\begin{tabular}{|c|}
\hline
$[m,\overline{p},f,h]$\\
\hline 
\hline 
$[ 20, 9, 2, 4 ]$\\
$[ 21, 4, 3, 4 ]$\\
$[ 24, 11, 2, 4 ]$\\
$[ 24, 17, 2, 4 ]$\\
$[ 24, 19, 2, 4 ]$\\
$[ 28, 9, 3, 4 ]$\\
$[ 30, 19, 2, 4 ]$\\
$[ 33, 4, 5, 4 ]$\\
$[ 39, 4, 6, 4 ]$\\
$[ 39, 16, 3, 8 ]$\\
$[ 40, 13, 4, 4 ]$\\
$[ 40, 17, 4, 4 ]$\\
$[ 48, 11, 4, 4 ]$\\
$[ 48, 19, 4, 4 ]$\\
$[ 52, 17, 6, 4 ]$\\
$[ 55, 4, 10, 4 ]$\\
$[ 55, 16, 5, 8 ]$\\
$[ 56, 11, 6, 4 ]$\\
$[ 56, 23, 6, 4 ]$\\
$[ 56, 37, 6, 4 ]$\\
$[ 57, 4, 9, 4 ]$\\
$[ 60, 7, 4, 4 ]$\\
$[ 60, 13, 4, 4 ]$\\
$[ 60, 49, 2, 8 ]$\\
$[ 66, 25, 5, 4 ]$\\
$[ 68, 9, 8, 4 ]$\\
$[ 69, 4, 11, 4 ]$\\
$[ 70, 9, 6, 4 ]$\\
$[ 72, 11, 6, 4 ]$\\
$[ 72, 41, 6, 4 ]$\\
$[ 72, 43, 6, 4 ]$\\
$[ 77, 4, 15, 4 ]$\\
$[ 84, 11, 6, 4 ]$\\
$[ 84, 19, 6, 4 ]$\\
$[ 84, 53, 6, 4 ]$\\
$[ 84, 61, 6, 4 ]$\\
$[ 88, 3, 10, 4 ]$\\
$[ 88, 17, 10, 4 ]$\\
$[ 88, 19, 10, 4 ]$\\
$[ 90, 49, 6, 4 ]$\\
$[ 92, 9, 11, 4 ]$\\
$[ 93, 7, 15, 4 ]$\\
$[ 95, 4, 18, 4 ]$\\
$[ 96, 11, 8, 4 ]$\\
$[ 96, 19, 8, 4 ]$\\
$[ 99, 4, 15, 4 ]$\\
$[ 100, 9, 10, 4 ]$\\
$[ 104, 33, 12, 4 ]$\\
$[ 104, 37, 12, 4 ]$\\
$[ 105, 17, 12, 4 ]$\\
$[ 105, 37, 12, 4 ]$\\
$[ 111, 4, 18, 4 ]$\\
$[ 112, 11, 12, 4 ]$\\
$[ 112, 23, 6, 8 ]$\\
$[ 112, 37, 12, 4 ]$\\
$[ 114, 25, 9, 4 ]$\\
$[ 116, 5, 14, 4 ]$\\
$[ 120, 17, 4, 8 ]$\\
\hline
\end{tabular}
\begin{tabular}{|c|}
\hline
$[m,\overline{p},f,h]$\\
\hline
\hline 
$[ 120, 53, 4, 8 ]$\\
$[ 120, 73, 4, 8 ]$\\
$[ 124, 9, 15, 4 ]$\\
$[ 129, 10, 21, 4 ]$\\
$[ 132, 5, 10, 4 ]$\\
$[ 132, 31, 10, 4 ]$\\
$[ 132, 47, 10, 4 ]$\\
$[ 136, 3, 16, 4 ]$\\
$[ 136, 5, 16, 4 ]$\\
$[ 138, 13, 11, 4 ]$\\
$[ 140, 3, 12, 4 ]$\\
$[ 140, 9, 6, 8 ]$\\
$[ 140, 17, 12, 4 ]$\\
$[ 141, 4, 23, 4 ]$\\
$[ 144, 11, 12, 4 ]$\\
$[ 144, 43, 12, 4 ]$\\
$[ 147, 4, 21, 4 ]$\\
$[ 148, 21, 18, 4 ]$\\
$[ 150, 19, 10, 4 ]$\\
$[ 152, 3, 18, 4 ]$\\
$[ 152, 33, 18, 4 ]$\\
$[ 152, 35, 18, 4 ]$\\
$[ 154, 9, 15, 4 ]$\\
$[ 155, 9, 30, 4 ]$\\
$[ 155, 41, 15, 8 ]$\\
$[ 156, 7, 12, 4 ]$\\
$[ 156, 11, 12, 4 ]$\\
$[ 156, 29, 6, 8 ]$\\
$[ 156, 49, 6, 8 ]$\\
$[ 161, 2, 33, 4 ]$\\
$[ 164, 5, 20, 4 ]$\\
$[ 165, 26, 10, 8 ]$\\
$[ 165, 37, 20, 4 ]$\\
$[ 165, 38, 20, 4 ]$\\
$[ 168, 11, 6, 8 ]$\\
$[ 168, 19, 6, 8 ]$\\
$[ 168, 59, 6, 8 ]$\\
$[ 168, 65, 6, 8 ]$\\
$[ 168, 67, 6, 8 ]$\\
$[ 168, 73, 6, 8 ]$\\
$[ 174, 13, 14, 4 ]$\\
$[ 176, 3, 20, 4 ]$\\
$[ 176, 19, 20, 4 ]$\\
$[ 177, 4, 29, 4 ]$\\
$[ 180, 7, 12, 4 ]$\\
$[ 180, 13, 12, 4 ]$\\
$[ 180, 49, 6, 8 ]$\\
$[ 183, 4, 30, 4 ]$\\
$[ 184, 3, 22, 4 ]$\\
$[ 184, 13, 22, 4 ]$\\
$[ 184, 31, 22, 4 ]$\\
$[ 188, 9, 23, 4 ]$\\
$[ 192, 11, 16, 4 ]$\\
$[ 192, 19, 16, 4 ]$\\
$[ 196, 9, 21, 4 ]$\\
$[ 198, 25, 15, 4 ]$\\
$[ 200, 13, 20, 4 ]$\\
$[ 200, 17, 20, 4 ]$\\
\hline
\end{tabular}
\begin{tabular}{|c|}
\hline
$[m,\overline{p},f,h]$\\
\hline
\hline 
$[ 201, 4, 33, 4 ]$\\
$[ 203, 4, 42, 4 ]$\\
$[ 203, 16, 21, 8 ]$\\
$[ 204, 5, 16, 4 ]$\\
$[ 204, 7, 16, 4 ]$\\
$[ 204, 11, 16, 4 ]$\\
$[ 204, 37, 16, 4 ]$\\
$[ 207, 4, 33, 4 ]$\\
$[ 209, 4, 45, 4 ]$\\
$[ 212, 9, 26, 4 ]$\\
$[ 213, 4, 35, 4 ]$\\
$[ 216, 11, 18, 4 ]$\\
$[ 216, 41, 18, 4 ]$\\
$[ 216, 43, 18, 4 ]$\\
$[ 219, 19, 36, 4 ]$\\
$[ 220, 3, 20, 4 ]$\\
$[ 220, 7, 20, 4 ]$\\
$[ 220, 9, 10, 8 ]$\\
$[ 220, 41, 10, 8 ]$\\
$[ 224, 11, 24, 4 ]$\\
$[ 224, 23, 12, 8 ]$\\
$[ 224, 37, 24, 4 ]$\\
$[ 228, 13, 18, 4 ]$\\
$[ 228, 43, 18, 4 ]$\\
$[ 228, 67, 18, 4 ]$\\
$[ 230, 9, 22, 4 ]$\\
$[ 231, 2, 30, 4 ]$\\
$[ 231, 4, 15, 8 ]$\\
$[ 231, 5, 30, 4 ]$\\
$[ 231, 19, 30, 4 ]$\\
$[ 232, 21, 28, 4 ]$\\
$[ 232, 73, 28, 4 ]$\\
$[ 237, 4, 39, 4 ]$\\
$[ 244, 5, 30, 4 ]$\\
$[ 248, 7, 30, 4 ]$\\
$[ 248, 19, 30, 4 ]$\\
$[ 248, 45, 30, 4 ]$\\
$[ 249, 4, 41, 4 ]$\\
$[ 253, 3, 55, 4 ]$\\
$[ 258, 13, 21, 4 ]$\\
$[ 264, 59, 10, 8 ]$\\
$[ 270, 49, 18, 4 ]$\\
$[ 276, 29, 22, 4 ]$\\
$[ 276, 31, 22, 4 ]$\\
$[ 276, 35, 22, 4 ]$\\
$[ 280, 37, 12, 8 ]$\\
$[ 280, 117, 12, 8 ]$\\
$[ 280, 137, 12, 8 ]$\\
$[ 282, 7, 23, 4 ]$\\
$[ 284, 9, 35, 4 ]$\\
$[ 286, 49, 30, 4 ]$\\
$[ 288, 11, 24, 4 ]$\\
$[ 288, 43, 24, 4 ]$\\
$[ 291, 25, 48, 4 ]$\\
$[ 292, 25, 36, 4 ]$\\
$[ 295, 4, 58, 4 ]$\\
$[ 296, 5, 36, 4 ]$\\
$[ 296, 17, 36, 4 ]$\\
\hline
\end{tabular}
\begin{tabular}{|c|}
\hline
$[m,\overline{p},f,h]$\\
\hline
\hline 
$[ 297, 4, 45, 4 ]$\\
$[ 299, 4, 66, 4 ]$\\
$[ 300, 13, 20, 4 ]$\\
$[ 300, 67, 20, 4 ]$\\
$[ 300, 109, 10, 8 ]$\\
$[ 304, 3, 36, 4 ]$\\
$[ 304, 35, 36, 4 ]$\\
$[ 308, 39, 30, 4 ]$\\
$[ 308, 135, 30, 4 ]$\\
$[ 308, 149, 30, 4 ]$\\
$[ 309, 4, 51, 4 ]$\\
$[ 310, 9, 30, 4 ]$\\
$[ 310, 41, 15, 8 ]$\\
$[ 312, 37, 12, 8 ]$\\
$[ 312, 41, 12, 8 ]$\\
$[ 312, 97, 12, 8 ]$\\
$[ 312, 149, 12, 8 ]$\\
$[ 316, 5, 39, 4 ]$\\
$[ 318, 7, 26, 4 ]$\\
$[ 321, 4, 53, 4 ]$\\
$[ 323, 9, 72, 4 ]$\\
$[ 327, 7, 27, 8 ]$\\
$[ 327, 28, 54, 4 ]$\\
$[ 328, 11, 40, 4 ]$\\
$[ 328, 13, 40, 4 ]$\\
$[ 329, 2, 69, 4 ]$\\
$[ 330, 7, 20, 4 ]$\\
$[ 330, 37, 20, 4 ]$\\
$[ 330, 49, 10, 8 ]$\\
$[ 330, 61, 10, 8 ]$\\
$[ 336, 11, 12, 8 ]$\\
$[ 336, 19, 12, 8 ]$\\
$[ 336, 59, 12, 8 ]$\\
$[ 336, 67, 12, 8 ]$\\
$[ 344, 3, 42, 4 ]$\\
$[ 344, 33, 42, 4 ]$\\
$[ 344, 67, 42, 4 ]$\\
$[ 345, 13, 44, 4 ]$\\
$[ 345, 17, 44, 4 ]$\\
$[ 348, 13, 14, 8 ]$\\
$[ 348, 19, 28, 4 ]$\\
$[ 348, 37, 28, 4 ]$\\
$[ 350, 9, 30, 4 ]$\\
$[ 352, 3, 40, 4 ]$\\
$[ 352, 19, 40, 4 ]$\\
$[ 354, 7, 29, 4 ]$\\
$[ 355, 4, 70, 4 ]$\\
$[ 355, 6, 35, 8 ]$\\
$[ 356, 5, 44, 4 ]$\\
$[ 357, 5, 48, 4 ]$\\
$[ 357, 37, 48, 4 ]$\\
$[ 360, 77, 12, 8 ]$\\
$[ 360, 97, 12, 8 ]$\\
$[ 360, 113, 12, 8 ]$\\
$[ 363, 4, 55, 4 ]$\\
$[ 368, 3, 44, 4 ]$\\
$[ 368, 13, 44, 4 ]$\\
$[ 368, 39, 22, 8 ]$\\
\hline
\end{tabular}
\begin{tabular}{|c|}
\hline
$[m,\overline{p},f,h]$\\
\hline
\hline
$[ 371, 4, 78, 4 ]$\\
$[ 372, 13, 30, 4 ]$\\
$[ 372, 41, 30, 4 ]$\\
$[ 372, 43, 30, 4 ]$\\
$[ 372, 59, 30, 4 ]$\\
$[ 376, 3, 46, 4 ]$\\
$[ 376, 7, 46, 4 ]$\\
$[ 376, 21, 46, 4 ]$\\
$[ 380, 3, 36, 4 ]$\\
$[ 380, 9, 18, 8 ]$\\
$[ 380, 23, 36, 4 ]$\\
$[ 381, 13, 63, 4 ]$\\
$[ 384, 11, 32, 4 ]$\\
$[ 384, 19, 32, 4 ]$\\
$[ 388, 25, 48, 4 ]$\\
$[ 390, 7, 12, 8 ]$\\
$[ 390, 107, 12, 8 ]$\\
$[ 392, 11, 42, 4 ]$\\
$[ 392, 23, 42, 4 ]$\\
$[ 392, 37, 42, 4 ]$\\
$[ 393, 4, 65, 4 ]$\\
$[ 395, 4, 78, 4 ]$\\
$[ 396, 5, 30, 4 ]$\\
$[ 396, 31, 30, 4 ]$\\
$[ 396, 47, 30, 4 ]$\\
$[ 402, 19, 33, 4 ]$\\
$[ 404, 9, 50, 4 ]$\\
$[ 406, 9, 42, 4 ]$\\
$[ 406, 23, 21, 8 ]$\\
$[ 407, 3, 90, 4 ]$\\
$[ 412, 17, 51, 4 ]$\\
$[ 413, 4, 87, 4 ]$\\
$[ 414, 13, 33, 4 ]$\\
$[ 417, 4, 69, 4 ]$\\
$[ 418, 5, 45, 4 ]$\\
$[ 420, 47, 12, 8 ]$\\
$[ 420, 53, 12, 8 ]$\\
$[ 420, 67, 12, 8 ]$\\
$[ 420, 73, 12, 8 ]$\\
$[ 420, 109, 6, 16 ]$\\
$[ 423, 4, 69, 4 ]$\\
$[ 424, 5, 52, 4 ]$\\
$[ 424, 33, 52, 4 ]$\\
$[ 426, 19, 35, 4 ]$\\
$[ 430, 9, 42, 4 ]$\\
$[ 432, 11, 36, 4 ]$\\
$[ 432, 43, 36, 4 ]$\\
$[ 436, 29, 54, 4 ]$\\
$[ 437, 4, 99, 4 ]$\\
$[ 438, 19, 36, 4 ]$\\
$[ 440, 13, 20, 8 ]$\\
$[ 440, 17, 20, 8 ]$\\
$[ 440, 37, 20, 8 ]$\\
$[ 440, 97, 20, 8 ]$\\
$[ 444, 19, 36, 4 ]$\\
$[ 444, 25, 18, 8 ]$\\
$[ 444, 35, 36, 4 ]$\\
$[ 448, 11, 48, 4 ]$\\
\hline
\end{tabular}
\begin{tabular}{|c|}
\hline
$[m,\overline{p},f,h]$\\
\hline
\hline 
$[ 448, 23, 24, 8 ]$\\
$[ 448, 37, 48, 4 ]$\\
$[ 450, 79, 30, 4 ]$\\
$[ 452, 9, 56, 4 ]$\\
$[ 453, 10, 75, 4 ]$\\
$[ 456, 67, 18, 8 ]$\\
$[ 460, 7, 44, 4 ]$\\
$[ 460, 9, 22, 8 ]$\\
$[ 460, 17, 44, 4 ]$\\
$[ 470, 9, 46, 4 ]$\\
$[ 471, 10, 78, 4 ]$\\
$[ 472, 3, 58, 4 ]$\\
$[ 472, 11, 58, 4 ]$\\
$[ 472, 33, 58, 4 ]$\\
$[ 473, 9, 105, 4 ]$\\
$[ 475, 4, 90, 4 ]$\\
$[ 483, 4, 33, 8 ]$\\
$[ 483, 10, 66, 4 ]$\\
$[ 483, 11, 66, 4 ]$\\
$[ 483, 26, 66, 4 ]$\\
$[ 488, 17, 60, 4 ]$\\
$[ 488, 157, 60, 4 ]$\\
$[ 489, 4, 81, 4 ]$\\
$[ 490, 9, 42, 4 ]$\\
$[ 492, 7, 40, 4 ]$\\
$[ 492, 11, 40, 4 ]$\\
$[ 492, 13, 40, 4 ]$\\
$[ 492, 17, 40, 4 ]$\\
$[ 495, 38, 60, 4 ]$\\
$[ 495, 58, 60, 4 ]$\\
$[ 495, 86, 30, 8 ]$\\
$[ 496, 19, 60, 4 ]$\\
$[ 496, 45, 60, 4 ]$\\
$[ 497, 2, 105, 4 ]$\\
$[ 498, 7, 41, 4 ]$\\
$[ 500, 9, 50, 4 ]$\\
$[ 501, 4, 83, 4 ]$\\
$[ 507, 4, 78, 4 ]$\\
$[ 507, 16, 39, 8 ]$\\
$[ 508, 9, 63, 4 ]$\\
$[ 510, 7, 16, 8 ]$\\
$[ 516, 19, 42, 4 ]$\\
$[ 516, 31, 42, 4 ]$\\
$[ 516, 61, 42, 4 ]$\\
$[ 517, 3, 115, 4 ]$\\
$[ 522, 13, 42, 4 ]$\\
$[ 525, 17, 60, 4 ]$\\
$[ 525, 37, 60, 4 ]$\\
$[ 528, 59, 20, 8 ]$\\
$[ 531, 4, 87, 4 ]$\\
$[ 536, 11, 66, 4 ]$\\
$[ 536, 19, 66, 4 ]$\\
$[ 536, 41, 66, 4 ]$\\
$[ 537, 4, 89, 4 ]$\\
$[ 539, 4, 105, 4 ]$\\
$[ 540, 7, 36, 4 ]$\\
$[ 540, 13, 36, 4 ]$\\
$[ 540, 49, 18, 8 ]$\\
\hline
\end{tabular}
\end{table}}
\newpage 
{\tiny 
\begin{table}[h]
\caption{The list of $(m,{\overline p},f)\in {\mathcal P}^\ast$ for $m\le 1000$ under the assumption that  $(h:=)|(\Z/m\Z)^\times/\langle p\rangle|> 2$ }
\label{Tab1}
\begin{tabular}{|c|}
\hline
$[m,\overline{p},f,h]$\\
\hline 
\hline 
$[ 543, 4, 90, 4 ]$\\
$[ 543, 13, 45, 8 ]$\\
$[ 548, 9, 68, 4 ]$\\
$[ 561, 29, 80, 4 ]$\\
$[ 561, 31, 80, 4 ]$\\
$[ 564, 7, 46, 4 ]$\\
$[ 564, 17, 46, 4 ]$\\
$[ 564, 59, 46, 4 ]$\\
$[ 568, 3, 70, 4 ]$\\
$[ 568, 15, 70, 4 ]$\\
$[ 568, 29, 70, 4 ]$\\
$[ 572, 15, 60, 4 ]$\\
$[ 572, 37, 60, 4 ]$\\
$[ 572, 49, 30, 8 ]$\\
$[ 573, 4, 95, 4 ]$\\
$[ 576, 11, 48, 4 ]$\\
$[ 576, 43, 48, 4 ]$\\
$[ 579, 25, 96, 4 ]$\\
$[ 581, 4, 123, 4 ]$\\
$[ 582, 25, 48, 4 ]$\\
$[ 582, 43, 24, 8 ]$\\
$[ 583, 4, 130, 4 ]$\\
$[ 584, 5, 72, 4 ]$\\
$[ 584, 11, 72, 4 ]$\\
$[ 588, 11, 42, 4 ]$\\
$[ 588, 53, 42, 4 ]$\\
$[ 588, 61, 42, 4 ]$\\
$[ 588, 103, 42, 4 ]$\\
$[ 594, 25, 45, 4 ]$\\
$[ 596, 9, 74, 4 ]$\\
$[ 597, 4, 99, 4 ]$\\
$[ 598, 49, 66, 4 ]$\\
$[ 600, 17, 20, 8 ]$\\
$[ 600, 53, 20, 8 ]$\\
$[ 600, 73, 20, 8 ]$\\
$[ 604, 5, 75, 4 ]$\\
$[ 605, 4, 110, 4 ]$\\
$[ 605, 16, 55, 8 ]$\\
$[ 606, 13, 50, 4 ]$\\
$[ 608, 3, 72, 4 ]$\\
$[ 608, 35, 72, 4 ]$\\
$[ 612, 5, 48, 4 ]$\\
$[ 612, 7, 48, 4 ]$\\
$[ 612, 11, 48, 4 ]$\\
$[ 612, 61, 48, 4 ]$\\
$[ 620, 7, 60, 4 ]$\\
$[ 620, 13, 60, 4 ]$\\
$[ 620, 21, 30, 8 ]$\\
$[ 620, 51, 30, 8 ]$\\
$[ 621, 4, 99, 4 ]$\\
$[ 627, 4, 45, 8 ]$\\
$[ 627, 13, 90, 4 ]$\\
$[ 627, 14, 90, 4 ]$\\
$[ 627, 17, 90, 4 ]$\\
$[ 628, 25, 78, 4 ]$\\
$[ 632, 5, 78, 4 ]$\\
$[ 632, 11, 78, 4 ]$\\
\hline
\end{tabular}
\begin{tabular}{|c|}
\hline
$[m,\overline{p},f,h]$\\
\hline 
\hline 
$[ 632, 31, 78, 4 ]$\\
$[ 633, 4, 105, 4 ]$\\
$[ 636, 19, 52, 4 ]$\\
$[ 636, 25, 26, 8 ]$\\
$[ 636, 61, 52, 4 ]$\\
$[ 638, 5, 70, 4 ]$\\
$[ 639, 4, 105, 4 ]$\\
$[ 642, 13, 53, 4 ]$\\
$[ 644, 3, 66, 4 ]$\\
$[ 644, 11, 66, 4 ]$\\
$[ 644, 37, 66, 4 ]$\\
$[ 644, 73, 66, 4 ]$\\
$[ 645, 13, 84, 4 ]$\\
$[ 645, 62, 84, 4 ]$\\
$[ 646, 9, 72, 4 ]$\\
$[ 646, 47, 36, 8 ]$\\
$[ 648, 11, 54, 4 ]$\\
$[ 648, 41, 54, 4 ]$\\
$[ 648, 43, 54, 4 ]$\\
$[ 649, 3, 145, 4 ]$\\
$[ 655, 4, 130, 4 ]$\\
$[ 655, 11, 65, 8 ]$\\
$[ 660, 17, 20, 8 ]$\\
$[ 660, 47, 20, 8 ]$\\
$[ 660, 53, 20, 8 ]$\\
$[ 660, 61, 10, 16 ]$\\
$[ 660, 83, 20, 8 ]$\\
$[ 663, 22, 48, 8 ]$\\
$[ 663, 29, 48, 8 ]$\\
$[ 664, 3, 82, 4 ]$\\
$[ 664, 19, 82, 4 ]$\\
$[ 664, 57, 82, 4 ]$\\
$[ 667, 4, 154, 4 ]$\\
$[ 667, 16, 77, 8 ]$\\
$[ 668, 9, 83, 4 ]$\\
$[ 669, 19, 111, 4 ]$\\
$[ 670, 19, 66, 4 ]$\\
$[ 672, 11, 24, 8 ]$\\
$[ 672, 19, 24, 8 ]$\\
$[ 672, 59, 24, 8 ]$\\
$[ 672, 67, 24, 8 ]$\\
$[ 676, 17, 78, 4 ]$\\
$[ 681, 4, 113, 4 ]$\\
$[ 687, 19, 57, 8 ]$\\
$[ 687, 46, 114, 4 ]$\\
$[ 688, 3, 84, 4 ]$\\
$[ 688, 67, 84, 4 ]$\\
$[ 690, 13, 44, 4 ]$\\
$[ 690, 17, 44, 4 ]$\\
$[ 692, 9, 86, 4 ]$\\
$[ 695, 4, 138, 4 ]$\\
$[ 696, 73, 28, 8 ]$\\
$[ 696, 77, 28, 8 ]$\\
$[ 696, 89, 28, 8 ]$\\
$[ 700, 3, 60, 4 ]$\\
$[ 700, 9, 30, 8 ]$\\
$[ 700, 17, 60, 4 ]$\\
\hline
\end{tabular}
\begin{tabular}{|c|}
\hline
$[m,\overline{p},f,h]$\\
\hline 
\hline 
$[ 704, 3, 80, 4 ]$\\
$[ 704, 19, 80, 4 ]$\\
$[ 705, 7, 92, 4 ]$\\
$[ 705, 23, 92, 4 ]$\\
$[ 708, 5, 58, 4 ]$\\
$[ 708, 7, 58, 4 ]$\\
$[ 708, 35, 58, 4 ]$\\
$[ 710, 9, 70, 4 ]$\\
$[ 710, 81, 35, 8 ]$\\
$[ 712, 3, 88, 4 ]$\\
$[ 712, 13, 88, 4 ]$\\
$[ 713, 9, 165, 4 ]$\\
$[ 714, 5, 48, 4 ]$\\
$[ 714, 11, 48, 4 ]$\\
$[ 714, 25, 24, 8 ]$\\
$[ 717, 4, 119, 4 ]$\\
$[ 723, 49, 120, 4 ]$\\
$[ 724, 33, 90, 4 ]$\\
$[ 726, 25, 55, 4 ]$\\
$[ 731, 9, 168, 4 ]$\\
$[ 732, 7, 60, 4 ]$\\
$[ 732, 35, 60, 4 ]$\\
$[ 732, 49, 30, 8 ]$\\
$[ 735, 17, 84, 4 ]$\\
$[ 735, 37, 84, 4 ]$\\
$[ 736, 3, 88, 4 ]$\\
$[ 736, 13, 88, 4 ]$\\
$[ 736, 39, 44, 8 ]$\\
$[ 737, 4, 165, 4 ]$\\
$[ 742, 9, 78, 4 ]$\\
$[ 744, 11, 30, 8 ]$\\
$[ 744, 19, 30, 8 ]$\\
$[ 744, 41, 30, 8 ]$\\
$[ 744, 43, 30, 8 ]$\\
$[ 744, 59, 30, 8 ]$\\
$[ 744, 73, 30, 8 ]$\\
$[ 747, 4, 123, 4 ]$\\
$[ 748, 3, 80, 4 ]$\\
$[ 748, 5, 80, 4 ]$\\
$[ 748, 7, 80, 4 ]$\\
$[ 748, 29, 80, 4 ]$\\
$[ 749, 4, 159, 4 ]$\\
$[ 750, 19, 50, 4 ]$\\
$[ 752, 3, 92, 4 ]$\\
$[ 752, 21, 92, 4 ]$\\
$[ 753, 7, 125, 4 ]$\\
$[ 755, 11, 75, 8 ]$\\
$[ 755, 34, 150, 4 ]$\\
$[ 759, 2, 110, 4 ]$\\
$[ 759, 5, 110, 4 ]$\\
$[ 759, 13, 110, 4 ]$\\
$[ 759, 37, 110, 4 ]$\\
$[ 760, 13, 36, 8 ]$\\
$[ 760, 17, 36, 8 ]$\\
$[ 760, 33, 36, 8 ]$\\
$[ 760, 93, 36, 8 ]$\\
$[ 764, 9, 95, 4 ]$\\
\hline
\end{tabular}
\begin{tabular}{|c|}
\hline
$[m,\overline{p},f,h]$\\
\hline 
\hline 
$[ 768, 11, 64, 4 ]$\\
$[ 768, 19, 64, 4 ]$\\
$[ 772, 25, 96, 4 ]$\\
$[ 776, 5, 96, 4 ]$\\
$[ 776, 59, 96, 4 ]$\\
$[ 784, 11, 84, 4 ]$\\
$[ 784, 23, 42, 8 ]$\\
$[ 784, 37, 84, 4 ]$\\
$[ 786, 7, 65, 4 ]$\\
$[ 788, 9, 98, 4 ]$\\
$[ 789, 4, 131, 4 ]$\\
$[ 790, 9, 78, 4 ]$\\
$[ 791, 9, 168, 4 ]$\\
$[ 792, 59, 30, 8 ]$\\
$[ 796, 9, 99, 4 ]$\\
$[ 799, 2, 184, 4 ]$\\
$[ 804, 7, 66, 4 ]$\\
$[ 804, 13, 66, 4 ]$\\
$[ 804, 19, 66, 4 ]$\\
$[ 805, 2, 132, 4 ]$\\
$[ 805, 17, 132, 4 ]$\\
$[ 808, 29, 100, 4 ]$\\
$[ 808, 73, 100, 4 ]$\\
$[ 810, 49, 54, 4 ]$\\
$[ 812, 11, 84, 4 ]$\\
$[ 812, 23, 42, 8 ]$\\
$[ 812, 45, 42, 8 ]$\\
$[ 812, 61, 84, 4 ]$\\
$[ 813, 4, 135, 4 ]$\\
$[ 820, 11, 40, 8 ]$\\
$[ 824, 7, 102, 4 ]$\\
$[ 824, 19, 102, 4 ]$\\
$[ 824, 29, 102, 4 ]$\\
$[ 828, 29, 66, 4 ]$\\
$[ 828, 31, 66, 4 ]$\\
$[ 828, 59, 66, 4 ]$\\
$[ 830, 9, 82, 4 ]$\\
$[ 831, 7, 138, 4 ]$\\
$[ 831, 10, 69, 8 ]$\\
$[ 834, 7, 69, 4 ]$\\
$[ 836, 17, 90, 4 ]$\\
$[ 836, 35, 90, 4 ]$\\
$[ 836, 47, 90, 4 ]$\\
$[ 840, 173, 12, 16 ]$\\
$[ 840, 193, 12, 16 ]$\\
$[ 846, 7, 69, 4 ]$\\
$[ 847, 4, 165, 4 ]$\\
$[ 849, 7, 141, 4 ]$\\
$[ 852, 19, 70, 4 ]$\\
$[ 852, 29, 70, 4 ]$\\
$[ 852, 83, 70, 4 ]$\\
$[ 856, 3, 106, 4 ]$\\
$[ 856, 17, 106, 4 ]$\\
$[ 856, 43, 106, 4 ]$\\
$[ 860, 9, 42, 8 ]$\\
$[ 860, 13, 84, 4 ]$\\
$[ 860, 23, 84, 4 ]$\\
\hline
\end{tabular}
\begin{tabular}{|c|}
\hline
$[m,\overline{p},f,h]$\\
\hline 
\hline 
$[ 861, 17, 120, 4 ]$\\
$[ 861, 58, 120, 4 ]$\\
$[ 864, 11, 72, 4 ]$\\
$[ 864, 43, 72, 4 ]$\\
$[ 869, 4, 195, 4 ]$\\
$[ 872, 13, 108, 4 ]$\\
$[ 872, 57, 108, 4 ]$\\
$[ 876, 5, 72, 4 ]$\\
$[ 876, 19, 36, 8 ]$\\
$[ 876, 31, 72, 4 ]$\\
$[ 876, 77, 18, 16 ]$\\
$[ 884, 15, 24, 16 ]$\\
$[ 885, 7, 116, 4 ]$\\
$[ 885, 17, 116, 4 ]$\\
$[ 888, 5, 36, 8 ]$\\
$[ 888, 13, 36, 8 ]$\\
$[ 888, 17, 36, 8 ]$\\
$[ 888, 217, 36, 8 ]$\\
$[ 891, 4, 135, 4 ]$\\
$[ 892, 9, 111, 4 ]$\\
$[ 893, 4, 207, 4 ]$\\
$[ 894, 7, 74, 4 ]$\\
$[ 895, 4, 178, 4 ]$\\
$[ 896, 11, 96, 4 ]$\\
$[ 896, 23, 48, 8 ]$\\
$[ 896, 37, 96, 4 ]$\\
$[ 897, 2, 132, 4 ]$\\
$[ 897, 4, 66, 8 ]$\\
$[ 897, 11, 132, 4 ]$\\
$[ 897, 61, 66, 8 ]$\\
$[ 900, 13, 60, 4 ]$\\
$[ 900, 67, 60, 4 ]$\\
$[ 900, 169, 30, 8 ]$\\
$[ 904, 3, 112, 4 ]$\\
$[ 904, 5, 112, 4 ]$\\
$[ 912, 67, 36, 8 ]$\\
$[ 913, 3, 205, 4 ]$\\
$[ 915, 97, 60, 8 ]$\\
$[ 915, 107, 60, 8 ]$\\
$[ 916, 5, 114, 4 ]$\\
$[ 917, 4, 195, 4 ]$\\
$[ 920, 13, 44, 8 ]$\\
$[ 920, 37, 44, 8 ]$\\
$[ 920, 73, 44, 8 ]$\\
$[ 921, 7, 153, 4 ]$\\
$[ 924, 5, 30, 8 ]$\\
$[ 924, 17, 30, 8 ]$\\
$[ 924, 25, 15, 16 ]$\\
$[ 924, 53, 30, 8 ]$\\
$[ 924, 61, 30, 8 ]$\\
$[ 924, 149, 30, 8 ]$\\
$[ 924, 157, 30, 8 ]$\\
$[ 924, 193, 30, 8 ]$\\
$[ 932, 9, 116, 4 ]$\\
$[ 933, 4, 155, 4 ]$\\
$[ 939, 13, 156, 4 ]$\\
$[ 940, 9, 46, 8 ]$\\
\hline
\end{tabular}
\begin{tabular}{|c|}
\hline
$[m,\overline{p},f,h]$\\
\hline 
\hline 
$[ 940, 13, 92, 4 ]$\\
$[ 940, 23, 92, 4 ]$\\
$[ 943, 2, 220, 4 ]$\\
$[ 944, 3, 116, 4 ]$\\
$[ 944, 11, 116, 4 ]$\\
$[ 946, 9, 105, 4 ]$\\
$[ 948, 5, 78, 4 ]$\\
$[ 948, 7, 78, 4 ]$\\
$[ 948, 11, 78, 4 ]$\\
$[ 948, 37, 78, 4 ]$\\
$[ 952, 3, 48, 8 ]$\\
$[ 952, 5, 48, 8 ]$\\
$[ 952, 11, 48, 8 ]$\\
$[ 952, 23, 48, 8 ]$\\
$[ 952, 31, 48, 8 ]$\\
$[ 952, 37, 48, 8 ]$\\
$[ 952, 65, 48, 8 ]$\\
$[ 952, 73, 48, 8 ]$\\
$[ 954, 7, 78, 4 ]$\\
$[ 955, 4, 190, 4 ]$\\
$[ 955, 16, 95, 8 ]$\\
$[ 956, 5, 119, 4 ]$\\
$[ 957, 2, 140, 4 ]$\\
$[ 957, 31, 140, 4 ]$\\
$[ 959, 2, 204, 4 ]$\\
$[ 963, 4, 159, 4 ]$\\
$[ 964, 29, 120, 4 ]$\\
$[ 966, 11, 66, 4 ]$\\
$[ 966, 19, 66, 4 ]$\\
$[ 966, 25, 33, 8 ]$\\
$[ 966, 59, 66, 4 ]$\\
$[ 968, 17, 110, 4 ]$\\
$[ 968, 19, 110, 4 ]$\\
$[ 968, 59, 110, 4 ]$\\
$[ 978, 43, 81, 4 ]$\\
$[ 979, 5, 220, 4 ]$\\
$[ 980, 3, 84, 4 ]$\\
$[ 980, 9, 42, 8 ]$\\
$[ 980, 17, 84, 4 ]$\\
$[ 987, 2, 138, 4 ]$\\
$[ 987, 10, 138, 4 ]$\\
$[ 987, 17, 138, 4 ]$\\
$[ 987, 58, 138, 4 ]$\\
$[ 989, 9, 231, 4 ]$\\
$[ 990, 7, 60, 4 ]$\\
$[ 990, 49, 30, 8 ]$\\
$[ 990, 61, 30, 8 ]$\\
$[ 990, 97, 60, 4 ]$\\
$[ 992, 19, 120, 4 ]$\\
$[ 992, 45, 120, 4 ]$\\
$[ 993, 19, 165, 4 ]$\\
$[ 995, 4, 198, 4 ]$\\
$[ 996, 7, 82, 4 ]$\\
$[ 996, 11, 82, 4 ]$\\
$[ 996, 17, 82, 4 ]$\\
$[ 1000, 13, 100, 4 ]$\\
$[ 1000, 17, 100, 4 ]$\\
\hline
\end{tabular}
\end{table}}

\end{document}